\documentclass{amsart}
\usepackage[dvips]{epsfig}
\usepackage{graphicx}
\usepackage{latexsym}
\usepackage{amsmath}
\usepackage{amsthm}
\usepackage{amssymb}
\RequirePackage[colorlinks=true,citecolor=blue,urlcolor=blue]{hyperref}
\usepackage{scrextend}
\usepackage{hyperref}

\usepackage{xcolor}

\usepackage{caption}
\usepackage{subcaption}
\usepackage{float}

 \usepackage[applemac]{inputenc}

\newtheorem{assumption}{Assumption}

\newtheorem{proposition}{Proposition}
\newtheorem{remark}{Remark}

\newcommand{\vip}{\vskip.2cm}

\newcommand{\COMMENTAIRE}[1]{}

\newcommand\numberthis{\addtocounter{equation}{1}\tag{\theequation}}

\begin{document}

\title[]{Local polynomial estimation of the intensity of a doubly stochastic Poisson process with bandwidth selection procedure}

\author{Thomas Deschatre}

\address{Thomas Deschatre, EDF Lab, Boulevard Gaspard Monge, Palaiseau, France.}

\email{thomas.deschatre@gmail.com}

\begin{abstract} 
We consider a doubly stochastic Poisson process with stochastic intensity $\lambda_t =n q\left(X_t\right)$ where $X$ is a continuous It\^o semimartingale and $n$ is an integer. Both processes are observed continuously over a fixed period $\left[0,T\right]$. An estimation procedure is proposed  in a non parametrical setting for the function $q$ on an interval $I$ where $X$ is sufficiently observed using a local polynomial estimator. A method to select the bandwidth in a non asymptotic framework is proposed, leading to an oracle inequality. If $m$ is the degree of the chosen polynomial, the accuracy of our estimator over the H\"older class of order $\beta$ is $n^{\frac{-\beta}{2\beta+1}}$ if $m \geq \lfloor \beta \rfloor$ and it is optimal in the minimax sense if $m \geq \lfloor \beta \rfloor$. A parametrical test is also proposed to test if $q$ belongs to some parametrical family. Those results are applied to French temperature and electricity spot prices data where we infer the intensity of electricity spot spikes as a function of the temperature.
\end{abstract}

\maketitle




\section{Introduction}

\subsection{Motivation} Jump processes and point processes are used in several domains such as finance, insurance or neuroscience, see \cite{karr91} for more details. In finance, they allow to model discontinuities in equity prices time series and heavy tails in asset returns \cite{tankov03}. An application for insurance is the model of extreme events such that occurrence times of earthquakes \cite{ogata86}. In neuroscience, these processes model spikes which is a potential difference in the membrane of a neuron \cite{truccolo05}. Spikes is a high increase of the potential in the membrane followed by a quick reversion to the initial level of the potential. They are also present in electricity spot prices time series and can be both negative and positive: price level is very high or very low during a short time period before coming back to its original level \cite{benth12}. One way to model them is to use mean reverting jump processes \cite{cartea05, meyer08}.

\medskip
In all these areas, the frequency of the jumps can be explained by an exogenous variable. Modeling these dependences can not be omitted, because they have an impact on risk management or prediction and help us to understand some behaviors. In \cite{truccolo05}, the author explains the neural spiking activity with three king of covariates: the previous spikes, exogenous stimuli and concurrent neural activity. In \cite{ait15}, the authors propose a model for financial contagion. Financial contagion is the fact that a large price move in a market causes large price moves in other markets ; jumps are explained by a first jump and in this case we often use multidimensional Hawkes process which are mutually exciting processes, see \cite{ait15}. For electricity spot prices, spikes are often caused by abnormal temperatures which are not modeled by a jump process. In a general case, when the covariate is not an other jump process, we often use doubly stochastic Poisson processes, which are Poisson processes with stochastic intensity. Two of the most famous models are the Aalen multiplicative model introduced in \cite{aalen78} where the intensity process is of the form $\alpha_t Y_t$ with $\alpha_t$ a function of time and $Y_t$ a stochastic process and the Cox regression model introduced in \cite{cox92} where the intensity process is of the form $\alpha_t \exp\left(\beta^T Z\right)$ with $Z$ a multi-dimensional random variable. A non parametric version of the Cox model exists where the intensity is of the form $\alpha_t \exp\left(f\left(Z\right)\right)$ \cite{castellan00}. These models are used mainly for life times modeling.

\medskip
A large literature is dedicated on methods of estimation for intensity estimation of Poisson process, especially in a non parametric setting. In the case of inhomogeneous Poisson processes, \cite{reynaud03} and \cite{reynaud14} use projection estimators and model selection techniques. Several finite dimensional spaces called models are considered to find projection estimators and they propose a penalty criterion in order to select a model. They work in a non asymptotic framework and a concentration inequality is derived. Furthermore, minimax rates are found over several classes of functions. A different approach for the estimation of the intensity is the use of kernel methods as in \cite{diggle85} and \cite{brooks91}. In \cite{diggle85} and \cite{brooks91} , asymptotic properties of the kernel estimator are studied ; in \cite{brooks91}, methods to select the bandwidth is proposed. In \cite{zhang10}, the intensity can be stochastic and is also estimated as a function of time with a kernel estimator in an asymptotic framework. In the context of Cox and Aalen processes, \cite{comte11} also proposes model selection techniques with projection estimators ; local polynomial estimator, which is a generalization of kernel estimators, is proposed by \cite{chen11} and studied in an asymptotic framework. A method of estimation in asymptotic framework for Cox regression with a time dependent covariable is established in \cite{murphy91}. Lasso penalization is proposed in \cite{lemler16} for Cox regression in the context of high dimensional covariate.

\subsection{Objectives and results} In our case, we are interested in a doubly stochastic Poisson process denoted by $N$ where its intensity $\lambda$ is a function of an exogenous covariate which is a stochastic process $X_t$ 
\[\lambda_t = q\left(X_t\right)\]
and our goal is to estimate the function $q$. In this case, conditionally on $\left(X_s\right)_{s \geq 0}$, $\left(N_s\right)_{s \geq 0}$ is a inhomogeneous Poisson process. We assume that we observe $N$ and $X$ over a time horizon $\left[0,T\right]$. We can think of the example of the frequency of electricity spot prices spikes as a function of the temperature depending also on time.

\medskip
This framework has already been studied. Indeed, \cite{utikal93} proposes a kernel estimator of the function $q$ in the case where $T$ goes to $\infty$ and when $X$ satisfies some asymptotical conditions, which can be for instance stationarity. A kernel estimator is also proposed in \cite{nielsen95} where the covariates depends on time and asymptotic properties are studied in the case of $n$ i.i.d. observations and $n$ goes to infinity. In the same context of i.i.d. observations with time dependent covariables, \cite{osullivan93} proposes a non parametric estimator based on model selection and asymptotic properties are studied. In \cite{delattre13}, $X_t$ corresponds to the fractional part of a Brownian motion and the doubly stochastic Poisson process is used to model the limit order book ; an estimation procedure is proposed in an asymptotic framework.

\medskip
We consider a different framework where $X$ is a continuous It\^o semimartingale having a local time $l_T^x$ for $x$ in $\mathbb{R}$. This local time measures the time spend by $X$ around $x$ before time $T$ and verifies properties of Proposition \ref{propertiesX}. It existence and properties can be insured by low restrictions given in Assumption \ref{assumptionsX}. The function $q$ can be estimated at point $x$ only if $X$ takes this value before time $T$, or if $l_T^x > 0$. We estimate $q$ on an arbitrarily interval $I$ and we work conditionally on the event 
\[D\left(I,\nu \right) = \{\omega \in \Omega, \; \underset{x \in I}{\inf} \;l^x_T\left(\omega\right) \geq \frac{\nu T}{|I|}\}\]
with $\nu \in \left(0,1\right]$. We choose to work in a non parametric framework and in a non asymptotic framework. To our knowledge, inferencing the intensity of a doubly stochastic Poisson process as a function of a continuous It\^o semimartingale in a non-asymptotic framework is not present in the literature. 
%
%

\medskip
A local polynomial estimator $\hat{q}_h$ of $q$ is proposed in Section \ref{kernelestimation} with $h$ a bandwidth parameter. The criteria used to evaluate the performance of our estimator is the $L_2$ norm on $I$, $\|\cdot\|_I$. We also gives a method to select a bandwidth over a finite set $\mathcal{H}$. We adapt the method of \cite{lacour16} used for density estimation with i.i.d. observations to our context of intensity estimation for doubly stochastic Poisson process. The method consists in approximating the bias of $\hat{q}_h$ by an estimator of $\|q_h - q_{h_{\min}}\|^2_I$ where $q_h = \mathbb{E}\left(\hat{q}_h | \mathcal{F}^X_T\right)$, $\left(\mathcal{F}^X_t\right)_{0 \leq t \leq T}$ is the natural filtration of $X$ and $h_{\min} = \min \; \mathcal{H}$. Indeed, if $h_{\min}$ is sufficiently small, the bias of $\hat{q}_{h_{\min}}$ is negligible and $\|q_h - q_{h_{\min}}\|^2_I \approx \|q-q_h\|^2_I$. A biased estimator of $\|q_h - q_{h_{\min}}\|^2_I$ is $\|\hat{q}_h - \hat{q}_{h_{\min}}\|^2_I$. Correcting this bias and adding an estimator of the variance term, our method consists in choosing a bandwidth $\hat{h}$ minimizing a criteria of the form 
\[\|\hat{q}_h - \hat{q}_{h_{\min}}\|^2_I + pen_{\alpha}\left(h\right)\]
where $pen_{\alpha}\left(h\right)$ is a penalty function and $\alpha > 0$ a parameter chosen by the statistician weighting the variance. This method is an extension of the one of Goldenshluger and Lepski \cite{goldenshluger11}.

\medskip
One of our main results is the oracle inequality of Proposition \ref{oracle}. If we write 
\[\lambda_t = n q\left(X_t\right)\] 
with $n \in \mathbb{N}, \; n \geq 1$, we obtain
\begin{align*}
\mathbb{E}\left(\|q-\hat{q}_{\hat{h}}\|^2_I | D\left(I,\nu\right)\right) &\leq \left(\alpha \vee \frac{1}{\alpha}+O\left(\log\left(n\right)^{-1}\right)\right) \underset{h \in \mathcal{H}}{\min}\; \mathbb{E}\left(\|q-\hat{q}_h\|^2_I | D\left(I,\nu\right)\right) \\
&+ O\left(\log\left(n\right)\mathbb{E}\left(\|q-q_{h_{\min}}\|^2_I | D\left(I,\nu\right)\right) \right)   + O\left(\frac{\log\left(n\vee |\mathcal{H}|\right)^6}{n}\right).
\end{align*}
This inequality is asymptotically optimal when $\alpha = 1$.

\medskip
Furthermore, if we consider that $q$ belongs the H\"older class on $I$ with parameter $\beta$, our estimator is optimal in the minimax rate sense and the minimax rate of convergence is $n^{-\frac{\beta}{2\beta+1}}$ if the degree of the polynomial is larger than $\lfloor \beta \rfloor$, see Proposition \ref{minimaxnonparametric}. In addition to give a method of estimation for the intensity as a function of $X$, we have also shown that the method of \cite{lacour16} is adapted to inhomogeneous Poisson process because the case $X_t = t$ respects the different assumptions.

\medskip 
A second objective is to test if our function $q$ belongs to some class of parametrical model. Indeed, non parametric estimators are not convenient for operational applications. This objective is achieved in Section \ref{parametricestimation} where a test is proposed. We test 
\[H_0: q = q_{\theta} \text{ on } I\; \text{for some }\theta \in \Theta\]
against 
\[H_1: q \neq q_{\theta} \text{ on } I\; \text{for all }\theta \in \Theta\]
where $\Theta \in \mathbb{R}^d$. We consider the contrast $M_n\left(\theta\right)$ defined in \eqref{contrast} that is an unbiased estimator of the distance between $q$ and $g_{\theta}$ based on our estimator $\hat{q}_{\hat{h}}$ which is an estimator of $q$ under both hypothesis. Under $H_0$, the minimum of the contrast gives an estimator of $\theta$, $\hat{\theta}_n$, converging at the rate $\sqrt{n}$ towards $\theta$ when $n \rightarrow \infty$ and a central limit theorem is provided, see Proposition \ref{parametricconvergence} (i). Furthermore, the quantity $\hat{h}^{\frac{1}{2}} n M_n\left(\hat{\theta}_n\right)$ converges in law towards a normal random variable under $H_0$, see Proposition \ref{parametricconvergence} (ii), but to $\infty$ under $H_1$, see Proposition \ref{parametricconvergence} (iii). This allows us to propose a critical region for the test.

\medskip
In Section \ref{application}, our estimation procedure is applied on electricity prices and temperature data in order to model the dependence between the spikes frequency of electricity prices and the temperature. In Section \ref{numericalresults}, results on simulated data are given.

\medskip
Proof of the oracle inequality is given in Section \ref{prooforacle} and other proofs in Section \ref{otherproofs}.

\section{Statistical setting}

Let $\left(X_t\right)_{0 \leq t \leq T}$ be a real valued continuous semimartingale of the form 
\begin{equation} \label{semimartingale}
X_t = X_0 + \int_0^t \mu_s ds + \int_0^t \sigma_s dW_s
\end{equation}
defined on a filtered probability space $\left(\Omega, \mathcal{F}, \left(\mathcal{F}_s\right)_{0 \leq s \leq T}, \mathbb{P}\right)$ where $\left(W_t\right)_{0 \leq t \leq T}$ is a standard Brownian motion, $\left(b_t\right)_{0 \leq t \leq T}$ and $\left(\sigma_t\right)_{0 \leq t \leq T}$ are c\'adl\'ag, progressively measurables and verify \\ 
$\int_0^T \left(|\mu_s| + \sigma^2_s\right) ds < \infty$ almost surely. Let also consider the natural filtration of $X$, $\left(\mathcal{F}^X_t\right)_{0 \leq t \leq T}$. Let $\left(N_t\right)_{0 \leq t \leq T}$ be a doubly stochastic Poisson process with intensity $\left(\lambda_t\right)_{0 \leq t \leq T}$ also defined on $\left(\Omega, \mathcal{F}, \left(\mathcal{F}_s\right)_{0 \leq s \leq T}, \mathbb{P}\right)$. We observe the two processes on $\left[0,T\right]$ with $T$ finite. We assume that the intensity, which is the function of interest, is of the form
\[\lambda_t =  n q\left(X_t\right)\]
where $n \in \mathbb{N}$, $n \geq 1$ corresponds to the asymptotic. As we observe the counting process on a finite time horizon, we need to have a sufficient number of jumps during this finite period, which is of order $\int_0^T \lambda_u du$ and then in our case of order $n$. We denote by $\Lambda_{\cdot} = \int_{0}^{\cdot} \lambda_u du$ the compensator of $N$ and by $M = N - \Lambda$ the compensated Poisson process. Conditionally on $\mathcal{F}^X_T$, $N$ is an inhomogeneous Poisson process with deterministic intensity at time $t$ $n q\left(X_t\right)$. Our aim is to estimate the function $\left(q\left(x\right), \; x \in I\right)$ for an arbitrary compact interval $I$.

\begin{assumption} \label{assumptionsX}
We assume that one of the following assumptions is true:
\begin{enumerate}
\item[(i)] $\underset{0 \leq s \leq T}{\inf} \sigma_s \geq \underline{\sigma}$ with $\underline{\sigma} > 0$ a deterministic constant and 
\[\mathbb{E}\left( \int_0^T |\mu_s| ds + \underset{0 \leq t \leq T}{\sup} |\int_0^t \sigma_s dW_s|\right) < \infty,\]
\item[(ii)] $X_t = t$ for all $t$ in $\left[0,T\right]$.
\end{enumerate}
\end{assumption} 

\begin{remark}
Assumption \ref{assumptionsX} (ii) can be more general. These assumptions are sufficient but not necessary conditions for existence of a local time. In the case $\mu$ deterministic and $\sigma = 0$, we need for all $0 \leq s \leq T$,  $\mu_s \neq 0$ but only the case $\mu = 1$ is interesting for us. Existence of local time with a stochastic drift and a null volatility could also be considered but existing results about local times for absolute continuous processes are not enough in the literature to consider it.   	
\end{remark}
 
\begin{proposition} \label{propertiesX} Let $X$ the process defined in \eqref{semimartingale}. Under Assumption \ref{assumptionsX}, there exists a function defined on $\mathbb{R} \times \left[0,T\right]$ and denoted $\left(x,t\right) \mapsto l_t^x$ verifying
\begin{itemize}
\item[(i)] an occupation time formula of the form
\[\int_0^t f\left(X_s\right)ds = \int_{\mathbb{R}} f\left(x\right) l^x_t dx, \; 0 \leq t \leq T\]
for any measurable function $f$ on $\Omega \times \mathbb{R}$, 
\item [(ii)]  $\mathbb{E}\left(\underset{x \in \mathbb{R}}{\sup}\; l^x_T \right) < \infty$ and
\item[(iii)] $x \mapsto l_T^x$ is continuous on $\mathbb{R}$ under Assumption \ref{assumptionsX} (i) and has one point of discontinuity under Assumption \ref{assumptionsX} (ii).
\end{itemize}
\end{proposition}

As noticed in \cite{hoffmann01}, the estimation of $q\left(x\right)$ at point $x \in I$ is meaningful only if the process $X$ hits the point $x$ before time $T$, or if $l^x_T > 0$. Indeed, $l_T^x$ is equal to 
\[\underset{\epsilon \to 0}{\lim} \; \frac{1}{2\epsilon} \int_{0}^T {\bf 1}_{|X_s-x| \leq \epsilon} ds\]
and measures the time spend by $X$ around the point $x$. For $\nu \in \left(0,1\right]$, let us define the event
\[D\left(I,\nu \right) = \{\omega \in \Omega, \; \underset{x \in I}{\inf} \;l^x_T\left(\omega\right) \geq \frac{\nu T}{|I|}\}\]
with $|I|$ the Lebesgue measure of $I$. From now, we work conditionally on the event $D\left(I,\nu\right)$, and assuming $\mathbb{P}\left(D\left(I,\nu\right)\right) > 0$. Under (ii), if $X_t = t$, the natural choice of $I$ is $\left[0,T\right]$ and $\nu = 1$.

\begin{remark}
We choose $\nu$ being dimensionless, justifying the normalization by $\frac{|I|}{T}$. Indeed, the local time have dimension equal to time times the inverse of the dimension of $X$. Furthermore, if $I = \left[\underline{I}, \bar{I}\right]$, and if we do the mapping $X'_t = \frac{X_{tT}-\underline{I}}{|I|}$ for $t \in \left[0,1\right]$, we have $\underset{x \in \left[0,1\right]}{\inf} \;l^x_1\left(X'\right) = \frac{|I|}{T} \underset{x \in I}{\inf} \;l^x_T$ where $l^x_1\left(X'\right)$ is the local time of $X'$ at time $1$ and point $x$. As $\int_I l_T^x dx \geq \underset{x \in I}{\inf} \;l^x_T |I|$ and $\int_I l_T^x dx = \int_{0}^T {\bf 1}_{X_s \in I} ds \leq T$, $\underset{x \in I}{\inf} \;l^x_T \leq \frac{T}{|I|}$ and $\nu$ has to be bounded by 1.
\end{remark}

\section{Local polynomial estimation}

\label{kernelestimation}

Let $m > 0$ be an integer, $K$ a kernel function and $\mathcal{S}^{+}_{m+1}$ the set of positive definite matrix of $\mathbb{R}^{(m+1) \times (m+1)}$. Let 
\[K_h\left(u\right) = h^{-1}K\left(\frac{u}{h}\right), \; u \in \mathbb{R}, \; h > 0.\]
Let us consider the local polynomial estimator, for $h > 0$ and $x \in \mathbb{R}$,
\begin{equation} \label{estimator}
\hat{q}_{h}\left(x\right) = \frac{1}{n} \int_{0}^T w\left(x,h,\frac{X_s-x}{h}\right) K_h\left(X_s-x\right){\bf 1}_{X_s \in I}  dN_s
\end{equation}
with 
\begin{equation} \label{w}
w\left(x,h,z\right) = U^T\left(0\right) B\left(x,h\right)^{-1} U\left(z\right) {\bf 1}_{B\left(x,h\right) \in \mathcal{S}^{+}_{m+1}}, \; z \in \mathbb{R},
\end{equation}
\[U\left(x\right) = \left(1,x,\frac{x^2}{2!},..,\frac{x^m}{m!}\right)^T,\]
and 
\begin{equation} \label{B}
B\left(x,h\right) =  \int_0^T U\left(\frac{X_s-x}{h}\right)  U^T\left(\frac{X_s-x}{h}\right) K_h\left(X_s-x\right){\bf 1}_{X_s \in I}  ds.
\end{equation}
If $B\left(x,h\right) \in \mathcal{S}^{+}_{m+1}$, the estimator $\hat{q}_h\left(x\right)$ is equal to $U^T\left(0\right) \hat{\theta}_h\left(x\right)$ with 
\[
\begin{split}
\hat{\theta}_h\left(x\right) &= \underset{\theta \in \mathbb{R}^{m+1}}{argmin} \;-\frac{2}{n}\theta^T \int_{0}^T U\left(\frac{X_s-x}{h}\right)K_h\left(X_s-x\right){\bf 1}_{X_s \in I}dN_s\\
&+ \theta^T\int_{0}^T U\left(\frac{X_s-x}{h}\right)U^T\left(\frac{X_s-x}{h}\right)K_h\left(X_s-x\right){\bf 1}_{X_s \in I} ds \theta.
\end{split}
\]
The case $l = 0$ corresponds to the classical Nadaraya Watson estimator. The term ${\bf 1}_{X_s \in I}$ allows us to avoid issues at the boundaries. We denote by $q_{h}$ the conditional expectation of $\hat{q}_{h}$ given $\mathcal{F}^X_T$: 
\[q_{h}\left(x\right) = \int_{0}^T w\left(x,h,\frac{X_s-x}{h}\right) K_h\left(X_s-x\right){\bf 1}_{X_s \in I}  q\left(X_s\right) ds.\]

\medskip
If $q$ is a polynomial function of degree $m$ on $I$, $q_h\left(x\right)$ is equal to $q\left(x\right)$, see Proposition \ref{polynom}.

\begin{proposition} \label{polynom}
Let $x \in I$ and $h > 0$ such that $B\left(x,h\right) \in \mathcal{S}^+_{m+1}$ where $B$ is defined in \eqref{B}. Let $Q$ be a polynomial function of degree $\leq m$. For any realization of the process $X$, we have 
\[ \int_{0}^T w\left(x,h,\frac{X_s-x}{h}\right) K_h\left(X_s-x\right){\bf 1}_{X_s \in I}  Q\left(X_s\right) ds = Q\left(x\right). 
\] 
In particular, 
\[
\int_{0}^T w\left(x,h,\frac{X_s-x}{h}\right) K_h\left(X_s-x\right){\bf 1}_{X_s \in I} ds = 1,
\]
\[ \int_{0}^T w\left(x,h,\frac{X_s-x}{h}\right) K_h\left(X_s-x\right){\bf 1}_{X_s \in I}\left(X_s-x\right)^k ds = 0, \; 1 \leq k \leq m.
\]
\end{proposition}

Proof of Proposition \ref{polynom} is immediate noticing that $\int_{0}^T w\left(x,h,\frac{X_s-x}{h}\right) K_h\left(X_s-x\right){\bf 1}_{X_s \in I}  Q\left(X_s\right) ds$ is the first component of the quantity minimizing 
\[\int_{0}^T \left(\theta^T U\left(\frac{X_s-x}{h}\right)-Q\left(X_s\right)\right)^2 K_h\left(X_s-x\right){\bf 1}_{X_s \in I}ds\]
and using the Taylor's expansion of a polynomial function ; it is similar to the proof of \cite[Proposition 1.12]{tsybakov09}.

\medskip 
For $B\left(x,h\right)^{-1}$ to be defined, the positive matrix $B\left(x,h\right)$ must be definite. Assumption \ref{assumptionkernel} is sufficient for $B\left(x,h\right)^{-1}$ to be well defined on the event $D\left(I,\nu\right)$, see Proposition \ref{minB}.

\begin{assumption} \label{assumptionkernel} We assume that:
\begin{enumerate}
\item[(i)] there exists $K_{\min} > 0$ and $\Delta > 0$ such that $K\left(u\right) \geq K_{\min} {\bf 1}_{|u| \leq \Delta} $ for all $u$ in $\mathbb{R}$,
\item[(ii)] $K$ has a compact support belonging to $\left[-1,1\right]$ and $\|K\|_{\infty} = \underset{x \in \mathbb{R}}{\sup} | K\left(x\right) | < \infty$.
\end{enumerate}
\end{assumption}

\begin{proposition} \label{minB} Let $0 < h \leq \frac{2}{3} \frac{|I|}{\Delta}$. Under Assumption \ref{assumptionkernel}, on the event $D\left(I,\nu\right)$, the matrix $B\left(x,h\right)$ defined in \eqref{B} belongs to $\mathcal{S}^+_{m+1}$ and for $x \in I$, $z \in \mathbb{R}$,
\[|w\left(x,h,z\right) {\bf 1}_{|z| \leq 1} | \leq \frac{|I|}{A_K \nu T}\]
with $w$ defined in Equation \eqref{w} and $A_{K}$ a constant depending on $K$.
\end{proposition}

\subsection{Method for bandwidth selection}

Our objective is to propose a method in order to choose the bandwidth $h$. We want for this bandwidth to minimize the $L_2$ loss on the interval $I$, $\mathbb{E}\left(\|q-\hat{q}_h\|^2_{I} | D\left(I,\nu\right)\right)$, with $\|f\|^2_I = \int_{I} f\left(x\right)^2dx$ for $f \in L_2\left(I\right)$ and let  $<\cdot,\cdot>_I$ the associated scalar product. This loss is equal to the sum of a bias term, $\mathbb{E}\left(\|q-q_h\|^2_{I} | D\left(I,\nu\right)\right)$, which depends on the regularity of $q$ and is usually increasing with $h$, and a variance term $\mathbb{E}\left(\|\hat{q}_h-q_h\|^2_{I} | D\left(I,\nu\right)\right)$, decreasing with $h$. The theoretical bandwidth minimizing this quantity depends on the function $q$ itself which is unknown and is called the oracle. One wants to find an estimator of this oracle, $\hat{h}$, such that $\mathbb{E}\left(\|q-\hat{q}_{\hat{h}}\|^2_{I} | D\left(I,\nu\right)\right) \leq \left(1+o\left(1\right)\right)  \underset{h \in \mathcal{H}}{\min} \mathbb{E}\left( \| \hat{q}_{h} - q\|_I^2| D\left(I,\nu\right)\right) + o\left(1\right)$ when $n \to \infty$, in order to have a loss with $\hat{h}$ close to the minimal one ; this type of inequality is called an oracle inequality.
The usual method to find this estimator of the oracle, which is done in this section, is to find an unbiased estimator of the bias and of the variance then to consider the bandwidth $h$ minimizing the sum of the two estimators.

\medskip 
In order to select the bandwidth parameter, we use the approach of \cite{lacour16} which is used for density estimation. We consider a finite set $\mathcal{H}$ of $\left(0,\infty\right)$ and $h_{\min} = \min \mathcal{H}$. The idea is to approximate the bias by $\mathbb{E}\left(\|q_h - q_{h_{\min}}\|_{I}^2 | D\left(I,\nu\right) \right)$ with $h_{\min}$ sufficiently small. A natural estimator of $\mathbb{E}\left(\|q_h - q_{h_{\min}}\|_{I}^2 | D\left(I,\nu\right)\right)$ is 
\begin{equation} \label{biasedestimatorbias}
\|\hat{q}_h - \hat{q}_{h_{\min}}\|^2_{I}.
\end{equation}
The estimator \eqref{biasedestimatorbias} induces a bias equal to the expectation of 
\[ \frac{1}{n} \int_0^T \int_{I} \left(w\left(x,h,\frac{X_s-x}{h}\right) K_h\left(X_s-x\right)-w\left(x,h_{\min},\frac{X_s-x}{h_{\min}}\right) K_{h_{\min}}\left(X_s-x\right)\right)^2 {\bf 1}_{X_s \in I} dx q\left(X_s\right)ds  \]
which can be estimated by the unbiased estimator
\[\frac{1}{n^2}\int_0^T \int_{I} \left(w\left(x,h,\frac{X_s-x}{h}\right) K_h\left(X_s-x\right)-w\left(x,h_{\min},\frac{X_s-x}{h_{\min}}\right) K_{h_{\min}}\left(X_s-x\right)\right)^2 {\bf 1}_{X_s \in I} dx dN_s.\]
This estimator can be written as 
\[\hat{V}_h + \hat{V}_{h_{\min}} - 2\hat{V}_{h,h_{\min}}\]
where 
\[\hat{V}_h = \frac{1}{n^2} \int_{0}^T \int_{I}  \left( w\left(x,h,\frac{X_s-x}{h}\right) K_h\left(X_s-x\right)\right)^2 {\bf 1}_{X_s \in I}dx dN_s\]
and 
\begin{align*} \label{Vhhmin}
&\hat{V}_{h,h_{\min}} = \\
&\frac{1}{n^2} \int_{0}^T \int_{I} w\left(x,h,\frac{X_s-x}{h}\right) K_h\left(X_s-x\right) w\left(x,h_{\min},\frac{X_s-x}{h_{\min}}\right) K_{h_{\min}}\left(X_s-x\right) {\bf 1}_{X_s \in I}dx dN_s\numberthis.
\end{align*}
An unbiased estimator of $\mathbb{E}\left(\|q_h - q_{h_{\min}}\|_{I}^2 | D\left(I,\nu\right)\right)$ is then 
\[\|\hat{q}_h - \hat{q}_{h_{\min}}\|^2_{I} -  \hat{V}_h - \hat{V}_{h_{\min}} + 2\hat{V}_{h,h_{\min}}.\]
An unbiased estimator of the variance which is equal to  
\[\mathbb{E}\left(\frac{1}{n} \int_{0}^T \int_{I}  \left(w\left(x,h,\frac{X_s-x}{h}\right) K_h\left(X_s-x\right)\right)^2 {\bf 1}_{X_s \in I} dx q\left(X_s\right)ds | D\left(I,\nu\right) \right)\]
is given by 
\[\hat{V}_h = \frac{1}{n^2} \int_{0}^T \int_{I}  \left( w\left(x,h,\frac{X_s-x}{h}\right) K_h\left(X_s-x\right)\right)^2 {\bf 1}_{X_s \in I}dx dN_s.\]
In order to choose the bandwidth, we then use the criteria 
\[\|\hat{q}_{h} - \hat{q}_{h_{\min}}\|_I^2 + pen_{\alpha}\left(h\right)\]
where 
\begin{equation} \label{pen}
pen_{\alpha}\left(h\right) = \alpha \hat{V}_h - \hat{V}_{h} - \hat{V}_{h_{\min}} + 2\hat{V}_{h,h_{\min}}, \text{ with } \alpha > 0.
\end{equation}
The term $\alpha > 0$ is used to weight the variance term. The optimal bandwidth $\hat{h}$ is given by 
\begin{equation} \label{optimalh}
\hat{h} = \underset{h \in \mathcal{H}}{\text{argmin }} \|\hat{q}_{h} - \hat{q}_{h_{\min}}\|_I^2 + pen_{\alpha}\left(h\right).
\end{equation}

In the following, we want to derive an oracle inequality for the estimator $\hat{q}_{\hat{h}}$ which has not be done to our knowledge.

\subsection{Concentration inequalities}

In order to compute this oracle inequality, we first need the two following concentration inequalities, from \cite{reynaud14} and \cite{houdre03}. The concentration inequality of Proposition \ref{bernstein} is a weak Bernstein inequality, the one of Proposition \ref{ustatistic} is an inequality for the Poisson U-statistic. These inequalities will be useful in our case because $N$ is an inhomogeneous Poisson process conditionally on $F^X_T$.

\begin{proposition}\cite[Equation (2.2)]{reynaud14} \label{bernstein} Let $T > 0$. Let $N$ be an inhomogeneous Poisson process with intensity $\lambda_{\cdot}$, $\Lambda_{\cdot} = \int_{0}^{\cdot} \lambda_u du$ and $M = N - \Lambda$. For all $u \geq 0$, with probability larger than $1-e^{-u}$
\[\int_{0}^T f\left(s\right) dM_s \leq \sqrt{2 u  \int_{0}^T f^2\left(s\right) d\Lambda_s} + \frac{\underset{x \in \left[0,T\right]}{\sup} |f\left(x\right)| u}{3}.\] 
\end{proposition}

\begin{proposition}\cite[Theorem 4.2]{houdre03} \label{ustatistic} Let $T > 0$. Let $N$ be an inhomogeneous Poisson process with intensity $\lambda_{\cdot}$, $\Lambda{\cdot} = \int_{0}^{\cdot} \lambda_u du$ and $M = N - \Lambda$. For all $\epsilon, u \geq 0$, with probability larger than $1-2.77e^{-u}$
\[\int_{0}^T \int_{0}^{s^-} f\left(u,s\right) dM_u dM_s \leq 2\left(1+\epsilon\right)^{\frac{3}{2}}C\sqrt{u} + 2\eta\left(\epsilon\right)Du + \beta\left(\epsilon\right)Bu^{\frac{3}{2}} + \gamma\left(\epsilon\right)Au^2\]
where 
\[\eta\left(\epsilon\right) = \sqrt{2\kappa}\left(2+\epsilon+\epsilon^{-1}\right), \; \beta\left(\epsilon\right) = e\left(1+\epsilon^{-1}\right)^2\kappa\left(\epsilon\right) + \left(\sqrt{2\kappa}\left(2+\epsilon+\epsilon^{-1}\right)\right)\vee \frac{\left(1+\epsilon\right)^2}{\sqrt{2}},\]
\[\gamma\left(\epsilon\right) = \left(e\left(1+\epsilon^{-1}\right)^2\kappa\left(\epsilon\right) \right)\vee \frac{\left(1+\epsilon\right)^2}{3},\; \kappa = 6, \;\kappa\left(\epsilon\right) = 1.25 + \frac{32}{\epsilon}\]
and 
\[A = \underset{\left(u,s\right) \in \left[0,T\right]^2}{\sup} f\left(u,s\right),\; B^2 = \max\{\underset{s \leq T}{\sup} \int_{0}^s f\left(u,s\right)^2 d\Lambda_u, \underset{u \leq T}{\sup} \int_{u}^T f\left(u,s\right)^2 d\Lambda_s   \},\]
\[C^2 = \int_{0}^T \int_{0}^s f\left(u,s\right)^2 d\Lambda_u d\Lambda_s, \; D = \underset{\int_{0}^T a^2_u d\Lambda_u = 1, \int_{0}^T b^2_s d\Lambda_s = 1}{\sup} \int_{0}^T a_u \int_{u}^T b_s f\left(u,s\right) d\Lambda_s d\Lambda_u.\]
\end{proposition}

\subsection{Oracle inequality}

For a function $f \in L_{\infty}\left(I\right)$, we denote by $\|f\|_{I,\infty}$ the norm $\underset{x \in I}{\sup} |f\left(x\right)|$. We will also need for the kernel the following norms: $\|\cdot\|_1$, $\|\cdot\| = \|\cdot\|_2$ and $\|\cdot\|_{\infty}$ corresponding respectively to the $L_1$, $L_2$ and $L_{\infty}$ norms on $\mathbb{R}$, with the $L_p$ norm defined by $\|f\|_p = \left(\int_{\mathbb{R}} |f\left(x\right)|^p \right)^{\frac{1}{p}}, \; p \geq 1$ and $\|f\|_{\infty} = \underset{x \in \mathbb{R}}{\sup} |f\left(x\right)|$. Proposition \ref{oracle} gives an oracle inequality for $\hat{q}_{\hat{h}}$.

\begin{proposition}\label{oracle} Assume \ref{assumptionsX} and \ref{assumptionkernel}. Let $x \geq 1$, $\epsilon \in \left(0,1\right)$. Let $\mathcal{H}$ a finite subset of $\left(0,\infty\right)$ such that $\min \mathcal{H} = h_{\min} \geq \frac{\|K\|_{\infty} \|K\|_1 |I|}{n}$ and $\max \mathcal{H} \leq  \frac{2}{3}\frac{|I|}{\Delta}$. Let $\hat{q}_h$ the local polynomial estimator defined in \eqref{estimator} and $\hat{h}$ the bandwidth defined in \eqref{optimalh}. With conditional probability given $\mathcal{F}^X_T$ larger than $C_1|\mathcal{H}|e^{-x}$, on the event $D\left(I,\nu\right)$, 
\begin{equation}  \label{oracleinequality}              
\begin{split}
\|\hat{q}_{\hat{h}}-q\|_I^2 &\leq C_0\left(\epsilon,\alpha\right) \underset{h \in \mathcal{H}}{\min}\; \| \hat{q}_{h} - q\|_I^2 +C_2\left(\epsilon,\alpha\right)\|q_{h_{\min}} - q\|_I^2 \\
&+ \frac{C_3\left(\epsilon,K,\alpha\right)|I|}{\nu^2 T^2}\left(\frac{\left(1+\|q\|_{I,\infty}\|l_T\|_{I,\infty}|I|\right)x^2}{n}  +\frac{x^3|I|}{n^2 h_{\min}}\right)  
\end{split}
\end{equation}
where $C_0\left(\epsilon,\alpha\right) = \alpha + \epsilon$ if $\alpha \geq 1$ and $C_0\left(\epsilon,\alpha\right) = \frac{1}{\alpha} + \epsilon$ if $0 < \alpha < 1$, $C_1$ is a constant, $C_2\left(\epsilon,\alpha\right)$ is a constant depending only on $\epsilon$ and $\alpha$ and $C_3\left(\epsilon,K,\alpha\right)$ is a constant depending only on $\epsilon$, $K$ and $\alpha$. Furthermore, $C_2\left(\epsilon,\alpha\right) \asymp \frac{1}{\epsilon}$ and $C_3\left(\epsilon,K,\alpha\right) \asymp \frac{1}{\epsilon^3}$ when $\epsilon \rightarrow 0$.

\medskip
We also have, 
\begin{equation} \label{oracleexpectation}
\begin{split}
\mathbb{E}(&\|\hat{q}_{\hat{h}}-q\|_I^2 | D\left(I,\nu\right)) \leq \left(\alpha \vee \frac{1}{\alpha} + \frac{\tilde{C}_1}{\log\left(n\right)}\right) \underset{h \in \mathcal{H}}{\min} \;\mathbb{E}\left( \| \hat{q}_{h} - q\|_I^2| D\left(I,\nu\right)\right)\\
& +\tilde{C}_2\left(\alpha\right)\log\left(n\right)\mathbb{E}\left(\|q_{h_{\min}} - q\|_I^2 | D\left(I,\nu\right)\right)\\
&+\frac{\tilde{C}_3\left(K,\alpha\right)|I|}{\nu^2 T^2}\left(\frac{\left(1+\|q\|_{I,\infty}\mathbb{E}\left(\|l_T\|_{I,\infty}| D\left(I,\nu\right)\right)|I|\right) \log\left(n\vee |\mathcal{H}|\right)^5}{n} +\frac{\log\left(n\vee |\mathcal{H}|\right)^6}{n}\right)  \\
&+\tilde{C}_4\frac{\|q\|^2_I}{n^4} + \frac{\tilde{C}_5\left(K\right)|I|}{\nu^2 T^2}\frac{\sqrt{\sum_{i=1}^4 \left(\|q\|_{\infty,I}T\right)^i}}{n}. 
\end{split}
\end{equation}
where $\tilde{C}_1$ and $\tilde{C}_4$ are constant, $\tilde{C}_2\left(\alpha\right)$ is a constant that only depends on $\alpha$, $\tilde{C}_3\left(K,\alpha\right)$ is a constant that only depends on $K$ and $\alpha$ and $\tilde{C}_5\left(K\right)$ is a constant that only depends on $K$. 
\end{proposition}

In inequality \eqref{oracleexpectation}, one can see the presence of an error of order $\log\left(n\right)\mathbb{E}\left(\|q_{h_{\min}} - q\|_I^2 | D\left(I,\nu\right)\right)$. This error comes from the approximation of the bias $\|q_h-q\|^2_I$ by $\|q_h-q_{h_{min}}\|^2_I$ and is negligible if $h_{\min}$ is small enough and $q$ regular enough. We also remark that the oracle inequality is asymptotically optimal when $\alpha = 1$.

\subsection{Adaptative minimax estimation}

In this section, we study the performance of the estimator $\hat{q}_{\hat{h}}$ in terms of convergence rate. We now work with the asymptotic $n \rightarrow \infty$, meaning the number of jumps becomes large when $n \rightarrow \infty$. For $\rho, \beta, L > 0$, let $\Lambda_{\rho,\beta} = \{f: I \rightarrow \mathbb{R}:\; f\left(x\right) \geq \rho, \; \|f\|_{I,\infty} < \infty \} \cap \Sigma\left(\beta,L,I\right)$ where $\Sigma\left(\beta,L,I\right)$ is the H\"older class on $I$ defined as the set of $l = \lfloor \beta \rfloor$ differentiable functions $f: I \to \mathbb{R}$ whose derivative $f^{\left(l\right)}$ verifies
\[|f^{\left(l\right)}\left(x\right)-f^{\left(l\right)}\left(x'\right)| \leq L |x-x'|^{\beta - l}, \; \forall x,\; x' \in I,\]
see \cite[Definition 1.2]{tsybakov09}. We will restrict to the study of $q \in \Lambda_{\rho,\beta}$. 

\medskip
To evaluate the performance of an estimator $\tilde{q}_n$ of $q$, we consider the minimax risk 
\[{\bf R}\left(\tilde{q}_n, \Lambda_{\rho,\beta}, \varphi_n\right) = \underset{q \in \Lambda_{\rho,\beta}}{\sup} \; \mathbb{E}\left( \varphi_n^{-2} \int_{I} \left(\tilde{q}_n\left(N,X,x\right)-q\left(x\right)\right)^2 dx | D\left(I,\nu\right)\right).\]
An estimator $\tilde{q}^*_n$ is said to attain an optimal rate of convergence $\varphi_n\left(\Lambda_{\rho,\beta}\right)$ if 
\[\underset{n \to \infty}{\lim \sup}\; {\bf R}\left(\tilde{q}^*_n, \Lambda_{\rho,\beta}, \varphi_n\left(\Lambda_{\rho,\beta}\right)\right) < \infty\]
and no estimator can attain a better rate:
\[\underset{n \to \infty}{\lim \inf} \; \underset{\tilde{q}_n}{\inf} \; {\bf R}\left(\tilde{q}_n, \Lambda_{\rho,\beta}, \varphi_n\left(\Lambda_{\rho,\beta}\right)\right) > 0\]
where the infimum is taken over all estimators.

\begin{proposition} \label{minimaxnonparametric} Assume \ref{assumptionsX} and \ref{assumptionkernel}. Let us consider the set of bandwidth $\mathcal{H} = \{h > 0 | h \geq \frac{\|K\|_{\infty}\|K\|_1 |I|}{n},\; h \leq \frac{2}{3}\frac{|I|}{\Delta} \text{ and } |I| h^{-1} \in \mathbb{N}\}$. Let $\hat{q}_{\hat{h}}$ the local polynomial estimator defined in Proposition \ref{oracle} and let $m$ be the degree of the corresponding polynomial. In the case where $m \geq \lfloor \beta \rfloor$, $\hat{q}_{\hat{h}}$ is optimal in the minimax sense and the optimal rate of convergence is given by $\varphi\left(\Lambda_{\rho,\beta}\right) = n^{\frac{-\beta}{2\beta+1}}$. In the case where $m < \lfloor \beta \rfloor$, the rate of convergence of $\hat{q}_{\hat{h}}$ is $n^{\frac{-m}{2m+1}}$
\end{proposition}

\section{Test for a parametric family}
\label{parametricestimation}
Let us consider the parametric family $P = \{g_{\theta}\left(\cdot\right), \; \theta \in \Theta\}$ with $\Theta$ a subset of $\mathbb{R}^d$, $d \geq 1$. Our objective is to test if the intensity function $q$ belongs to $P$. Let us consider the two hypothesis:
\[
\left\{
  \begin{array}{l}
H_0: \; \exists \theta_0 \in \Theta,\; q\left(\cdot\right) = g_{\theta}\left(\cdot\right) \\
H_1: \; q \notin P
 \end{array}
\right..
\]
We want to test $H_0$ against $H_1$. Under both hypothesis, one way to estimate $q$ is to use the local polynomial estimator $\hat{q}_{\hat{h}}$. As we work in an asymptotic framework, we denote by $h_n$ the optimal bandwidth $\hat{h}$. Under $H_0$, to estimate the parameter $\theta$, let us consider the contrast
\begin{equation}  \label{contrast}
\begin{split}
M_n\left(\theta\right) &= \| \hat{q}_{h_n}\left(\cdot\right) - \int_0^T w\left(\cdot,h_n,\frac{X_s-\cdot}{h_n}\right)K_{h_n}\left(X_s-\cdot\right) {\bf 1}_{X_s \in I} g_{\theta}\left(X_s\right)ds\|^2_I \\
&- \frac{1}{n^2} \int_I  \int_0^T w^2\left(x,h_n,\frac{X_s-x}{h_n}\right) K^2_{h_n}\left(X_s-x\right){\bf 1}_{X_s \in I}dN_s dx.
\end{split}
\end{equation}
The contrast defined in \eqref{contrast} is similar to the one in \cite{ait96}, used in the case of the estimation of the drift and the volatility of an It\^o diffusion. However, in \cite{ait96}, the norm is weighted by the density of $X$, that is $l_T^x$ is our case. As it is important to have a good estimate of $q$ everywhere on $I$, the norm is not weighted in \eqref{contrast}. The second term in the right hand side of \eqref{contrast} is a correction of the bias. We can also notice that we use the function $\int_0^T w\left(\cdot,h_n,\frac{X_s-\cdot}{h_n}\right)K_{h_n}\left(X_s-\cdot\right) {\bf 1}_{X_s \in I} g_{\theta}\left(X_s\right)ds$ and not directly $g_\theta\left(\cdot\right)$ in order to eliminate the bias term $\int_0^T w\left(\cdot,h_n,\frac{X_s-\cdot}{h_n}\right)K_{h_n}\left(X_s-\cdot\right) {\bf 1}_{X_s \in I} \left(g_{\theta}\left(X_s\right)-g_{\theta}\left(\cdot\right)\right)ds$ that would appear and then avoid us to make assumptions on the speed of convergence of this term which depends on the regularity of $g_{\theta}$. 
An estimator of $\theta$ under $H_0$ is  
\begin{equation} \label{estimatortheta}
\hat{\theta}_n = \underset{\theta \in \Theta}{\inf} M_n\left(\theta\right).
\end{equation}
Under classical Assumption \ref{assumpparametric}, this estimator is consistent at a speed rate of $\sqrt{n}$, see Proposition \ref{parametricconvergence} (i). The idea of the test is that under $H_0$, $M_n\left(\hat{\theta_n}\right)$ is close to $M\left(\theta_0\right)$ which is equal to 0. The rate of convergence is of order $n \sqrt{h_n}$, see Proposition \ref{parametricconvergence} (ii). However, under $H_1$, $M_n\left(\hat{\theta}_n\right)$ converges to $\underset{\theta \in \Theta}{\inf \;} \|q-g_{\theta}\|^2_I$ which is different from 0 and then $n \sqrt{h_n} M_n\left(\hat{\theta}_n\right)$ goes to $\infty$, see Proposition \ref{parametricconvergence} (iii). 

\begin{assumption} \label{assumpparametric}
We assume that 
\begin{enumerate}
\item[(i)] The set $\Theta$ is compact in $\mathbb{R}^d$.
\item[(ii)] For some $M = M_I > 0$,  
\[\underset{x \in I}{\sup} |g_{\theta_1}\left(x\right) - g_{\theta_2}\left(x\right)| \leq M \|\theta_1 - \theta_2\|_d \text{ for } \theta_1,\theta_2 \in \Theta\]
where $\|\cdot\|_d$ is the $\mathbb{R}^d$ Euclidian norm.
\item[(iii)] For all $x$ in $I$, $\theta \mapsto g_{\theta}\left(x\right)$ is three times continuously differentiable. Furthermore, $x \mapsto g_{\theta}\left(x\right)$, $x \mapsto \partial_{\theta} g_{\theta}\left(x\right)$ are continuous on $I$ and $\frac{\partial_{2,\theta} g_{\theta}}{\partial \theta_i \partial \theta_j}$, $\frac{\partial_{3,\theta} g_{\theta}}{\partial \theta_i \partial \theta_j \partial \theta_k}$ are bounded on $I$ for $i$, $j$, $k$ $=1,..,d$ for all $\theta \in \Theta$. 
\item[(iv)] For some $\eta = \eta_I > 0$, 
\[ \underset{\left(x,\theta\right) \in I \times \theta}{\inf} \lambda_{\min}\left( \partial_{\theta} g_\theta \left(x\right) \partial_{\theta} g_{\theta}\left(x\right)^T \right) \geq \eta\]
where $\lambda_{\min}\left(A\right)$ is the smallest eigenvalue of a matrix $A$.
\item[(v)] The equality $g_{\theta_1} = g_{\theta_2}$ on $I$ implies $\theta_1 = \theta_2$.
\item[(vi)] $h_n \rightarrow 0$ and $n \sqrt{h_n} \rightarrow \infty$.
\item[(vii)] The intensity function $q$ is continuous on $I$ and the kernel function $K$ is continuous on $\mathbb{R}$.
\end{enumerate}
\end{assumption}

\begin{proposition}  \label{parametricconvergence}
Let $M_n\left(\theta\right)$ and $\hat{\theta}_n$ defined respectively by \eqref{contrast} and \eqref{estimatortheta}. We work under Assumption \ref{assumptionsX}, Assumption \ref{assumptionkernel} and Assumption \ref{assumpparametric}. On the event $D\left(I,\nu\right)$, conditionally on $F^X_T$, under $H_0$, 
\begin{enumerate}
\item[(i)] 
\[\sqrt{n}\left(\hat{\theta}_n-\theta_0\right)  \overset{\mathcal{L}}{\rightarrow} \mathcal{N}\left(0,\left(\int_I \partial_{\theta} g_{\theta_0}\left(u\right)\partial_{\theta} g_{\theta_0}\left(u\right)^Tdu\right)^{-1} \int_I \frac{\partial_{\theta} g_{\theta_0}\left(u\right)\partial_{\theta} g_{\theta_0}\left(u\right)^T g_{\theta_0}\left(u\right)}{l_T^u}du\right),\]
\item[(ii)] 
\[ n \sqrt{h_n} M_n\left(\hat{\theta}_n\right) \overset{\mathcal{L}}{\rightarrow}  \mathcal{N}\left(0,2 \int_{\mathbb{R}} \left(\int_{\mathbb{R}} w\left(u\right)w\left(u+p\right)K\left(u\right) K\left(u+p\right) du\right)^2 dp \int_I \frac{\left(g_{\theta_0}\left(y\right)\right)^2}{\left(l_T^y\right)^2} dy\right)\]
\end{enumerate}
with
\[w\left(u\right) = U^T\left(0\right)\left(\int_{\mathbb{R}} U\left(z\right)U^T\left(z\right)K\left(z\right)dz\right)^{-1}U\left(u\right)\]
and under $H_1$, 
\begin{enumerate}
\item[(iii)] \[ |n \sqrt{h_n} M_n\left(\hat{\theta}_n\right)| \overset{p}{\rightarrow} \infty.\]
\end{enumerate}
\end{proposition}

In order to test the null hypothesis at level $\gamma \in \left(0,1\right)$, we reject $H_0$ when  
\[ |M_n\left(\theta_n\right)| \geq \hat{c}\left(\gamma\right) = n^{-1} h_n^{-\frac{1}{2}} \sqrt{\mathcal{V}_n} \Phi^{-1}\left(1-\frac{\gamma}{2}\right) \]
where 
\[
\mathcal{V}_n =A\left(K\right)\int_I \left(\frac{g_{\hat{\theta}_n}\left(y\right)}{\int_0^T K_{h_n}\left(y-X_s\right){\bf 1}_{X_s \in I} ds}\right)^2 dy,\]
 \[A\left(K\right) = 2\left(\int_{\mathbb{R}} K\left(u\right)du\right)^2 \int_{\mathbb{R}} \left(\int_{\mathbb{R}} w\left(u\right)w\left(u+p\right)K\left(u\right) K\left(u+p\right) du\right)^2 dp \]
only depends on $K$ and $\Phi$ is the cumulative distribution function of a $\mathcal{N}\left(0,1\right)$ random variable. 

\section{Dependence between the frequency of electricity spot spikes and temperature}
\label{application}
\subsection{Data}

We dispose of 
\begin{itemize}
\item the hourly French EPEX spot price between the first January of 2007 and the first of January 2017 not included,
\item the hourly French temperature, which is an spatial average of the temperature over 32 cities, between the first January of 2007 and the first of January 2017 not included.
\end{itemize}

\begin{figure}[h!]
             \centering
    \begin{subfigure}[b]{0.45\textwidth}
        \centering
        \includegraphics[width=\textwidth]{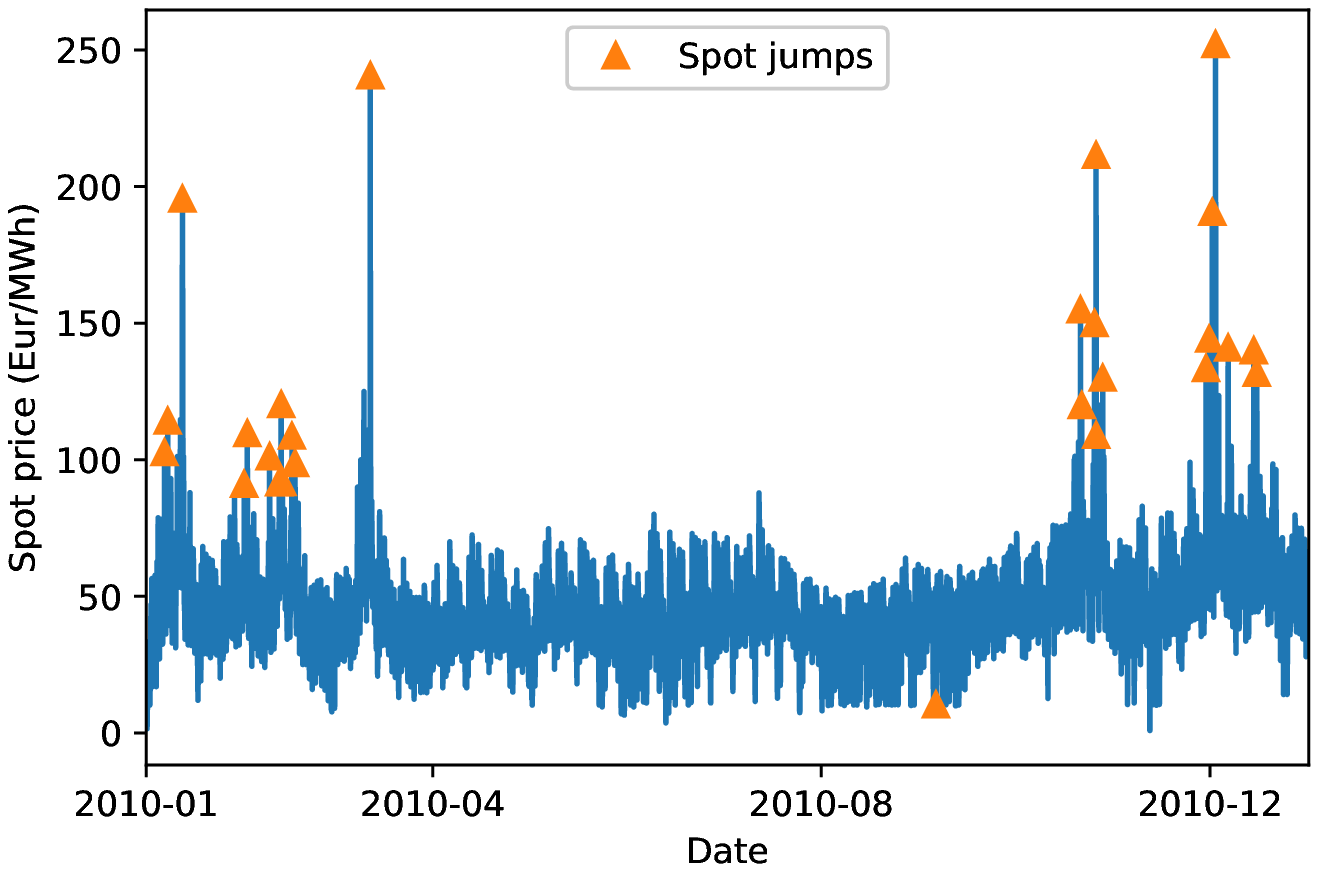}
        \caption{\it Spot price.}
    \end{subfigure}
    \begin{subfigure}[b]{0.45\textwidth}
        \centering
        \includegraphics[width=\textwidth]{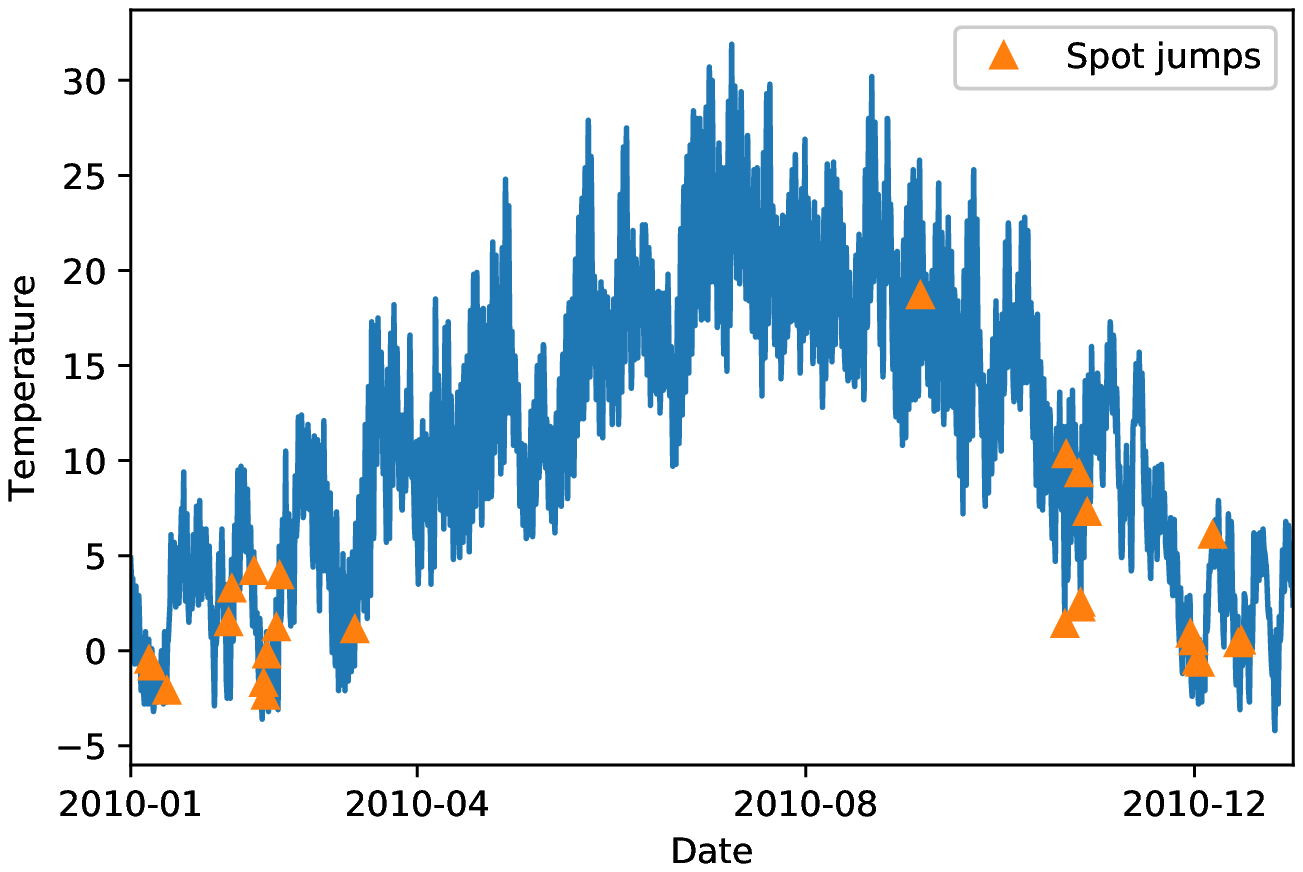}
        \caption{\it Temperature.}
    \end{subfigure}         
        \caption{\it \label{spottemperaturejumps} French spot price and temperature during 2010.}
    \end{figure}

The times series are given in Figure \ref{spottemperaturejumps} for the year 2010.

\subsection{Detection of the jumps} In the spot price time series, we observe spikes that are characteristic of the electricity spot market. A spike can be defined as a jump with a strong mean reversion. We then assume that the spot price $S$ has the following dynamic:
\[S_t = Y_t + Z_t\]
with $Y_t$ a continuous It\^o semi-martingale and 
\[Z_t = \int_0^t \int_{\mathbb{R}} x e^{-\beta\left(t-s\right)} \underline{p}\left(dt,dx\right)\]
with $\underline{p}$ a Poisson measure on $\mathbb{R}^+ \times \mathbb{R}$ with compensator $\underline{q} = \lambda_t dt \otimes \nu\left(dx\right)$. $Z$ has the dynamic 
\[dZ_t = -\beta Z_t dt + \int_{\mathbb{R}} x \underline{p}\left(dt,dx\right)\]
corresponding to a mean reverting Poisson process and models the spikes. Let us consider $N$ the Poisson process associated to the jump times of $X$. $N$ has intensity $\left(\lambda_t\right)_{0 \leq t \leq T}$ which is the intensity we want to estimate as a function of the temperature. In order to detect the jumps, we use the method of \cite{deschatre18} with a threshold equal to $5 \hat{\sigma} \Delta^{0.49}$ where $\Delta$ is the frequency of observations and $\hat{\sigma}$ is the multipower variation estimator of order 20. We keep only the increments verifying $\Delta_i S \Delta_{i+1} S < 0$ as jumps with $\Delta_i S = S_{t_i} - S_{t_{i-1}}$. When the frequency of observations $\Delta$ goes to 0 and when $\beta$ is large enough, \cite{deschatre18} proves that this filtering allows to detect with probability one every spikes under some asymptotic conditions. The data are segmented in periods of one year for the detection of the jumps in order to avoid too much change in the volatility. In \cite{deschatre18}, the intensity of the Poisson process is constant. Assuming that $\lambda$ is bounded below and above, the results can easily be extended to the case where $\lambda$ is stochastic. Jump times are represented in Figure \ref{spottemperaturejumps}. In the following, we consider that we observe $N$ but we are aware that we only have an estimator of it.

\subsection{Dependence with temperature}

In this section, we estimate the intensity of the jump process as a function of the temperature, using the method of Section \ref{kernelestimation}. In addition to the statistic interest, quanto options are financial options with temperature and spot price as underlying. They can be used for instance to hedge both volume and price risks. In order to price these options, it is necessary to capture the dependence between the temperature and the spot price. More details are given in \cite{benth15} where the dependence between the two is only modeled by a correlation and the spikes are not represented.

\begin{figure}[h!]
    \centering
         \includegraphics[width=0.7\textwidth]{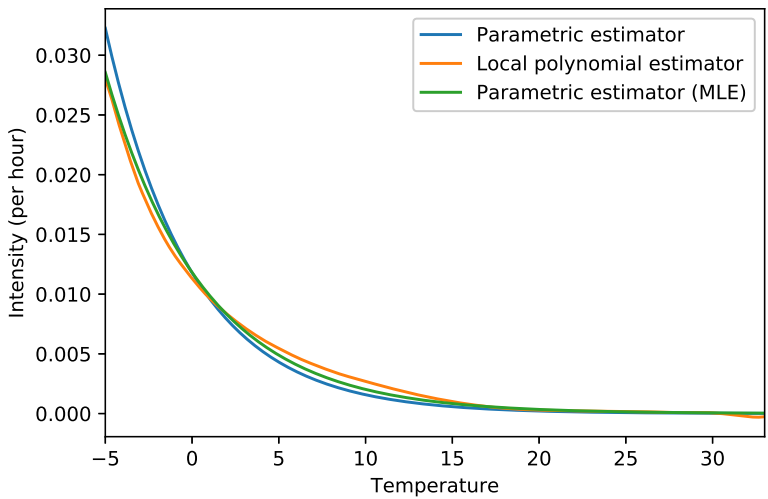}
        \caption{\label{estimatortemperature} \it Local polynomial estimator (in the case $\alpha \in \{0.5, 0.75, 1, 1.25, 1.5, 1.75\}$), parametric estimator and maximum likelihood parametric estimator (MLE) of the intensity as a function of the temperature.}
\end{figure}

\medskip
The temperature is illustrated in Figure \ref{spottemperaturejumps} along with the jump times. The spikes seems to happen more often for low temperatures. The observed temperature belongs to the interval $\left[-8.50,33.95\right]$. The temperature is not observed continuously but because of the high frequency of the data and the long range of observation, we pretend that the error due to the discretization is negligible. One wants to estimate the intensity of the spike process as a function of the temperature on the interval $\left[-5,33\right]$ where the temperature is sufficiently observed. To estimate the intensity function, we consider the Epanechnikov kernel $K\left(u\right) = \frac{3}{4}\left(1-u^2\right){\bf 1}_{|u| \leq 1}$ and the local polynomial estimator with degree 1 considered in Section \ref{kernelestimation}. We choose $h_{\min} $ equal to $\frac{|I| \|K\|_1 \|K\|_{\infty}}{N_I} = 0.13$, where  $\|K\|_1 = 1$, $\|K\|_{\infty} = \frac{3}{4}$ and $N_I = 219$ is the number of jumps in the interval $I$. The tuning parameter of the estimation procedure $\alpha$ is chosen equal to $1$. The optimal bandwidth is selected among the set $\mathcal{H} = \{h = h_{\min} + 0.1 i,\; h \leq 25 \}$. The minimum of the criteria is achieved for $\hat{h} = 8.73$ and the estimator for this value of $h$ is given in Figure \ref{estimatortemperature}. The estimator takes small negative values for high temperatures, which is caused by the total absence of jumps in this area ; one can take the maximum between the estimator and a small positive value for the intensity.

\begin{figure}[h!]
    \centering
         \includegraphics[width=0.7\textwidth]{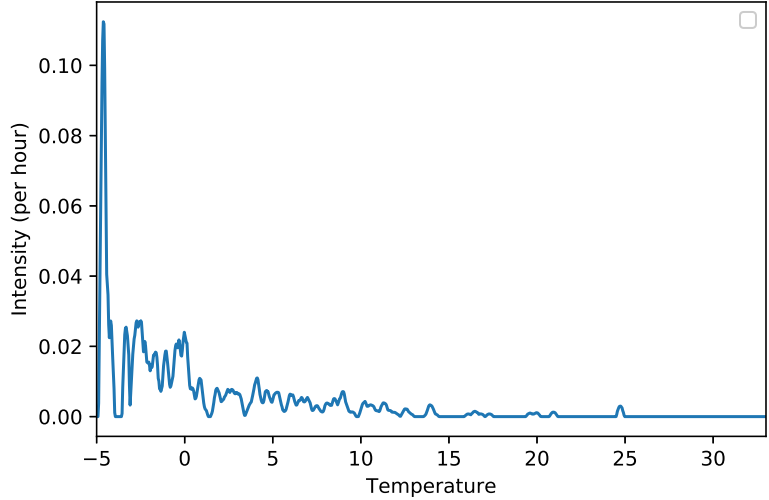}
        \caption{\label{estimatortemperature025} \it Local polynomial estimator of the intensity as a function of the temperature in the case $\alpha = 0.25$.}
\end{figure}   

\begin{figure}[h!]
    \centering
         \includegraphics[width=0.7\textwidth]{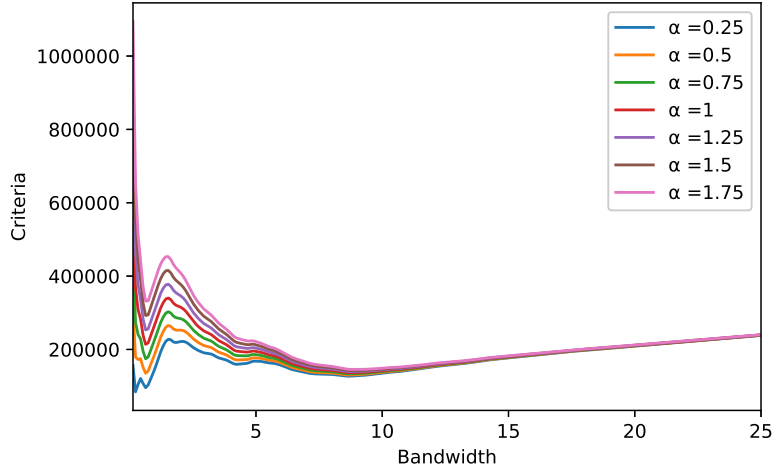}
        \caption{\label{criteriahyper} \it Criteria to minimize in order to find $\hat{h}$ for different values of $\alpha$.}
\end{figure}   

This result confirms our intuition: spikes happen more often when temperature is low. We now want to test the hypothesis that the intensity is a quadratic function of the temperature:
\[q\left(x\right) = a_0 \exp\left(a_1 x\right) \text{ for } x \in I\]
with $a_0 > 0$, $a_1 < 0$. The constant 
\[\int_{\mathbb{R}}\left(\int_{\mathbb{R}} w\left(u\right) w\left(u+p\right) K\left(u\right) K\left(u+p\right)du\right)^2 dp =  \frac{413113}{985600}\]
is needed for the test. We find that the null hypothesis is not rejected for a level of confidence at $95\%$ (with a p value equals to 0.083) and that the estimated parameters are $\left(\hat{a}_0,\hat{a}_1\right) = \left(1033.8,  -0.2\right)$. Figure \ref{estimatortemperature} includes the parametric estimator of the intensity as a function of the temperature in the case $\alpha = 1$. The maximum likelihood estimator (MLE) of the parametric model $q\left(x\right) = a_0 \exp\left(a_1 x\right)$ is equal to $\left(1035.7,  -0.17\right)$ and gives a similar estimator for $q$, see Figure \ref{estimatortemperature}. However, the MLE estimator of $q$ is closer to the local polynomial estimator than our parametric estimator ; this can be caused by the weight introduced in the norm in \eqref{contrast} for $g_{\theta}$. In a second time, we test if the intensity function is independent from the temperature, corresponding to $q$ constant: the test is rejected for a level of confidence at $95\%$ (with a p value equals to 0).  

\medskip
To study the sensibility to the choice of $\alpha$, we perform our estimation procedure for $\alpha \in \{0.25, 0.5, 0.75, 1.25, 1.5, 1.75\}$. Results remains the same for $\alpha \in \{0.5, 0.75, 1.25, 1.5, 1.75\}$: $\hat{h} = 8.73$. However, it differs significantly when $\alpha = 0.25$: $\hat{h} = 0.23$ which is closed to $h_{\min}$ leading to a high bias, which is consistent with a small value of $\alpha$, see Figure \ref{estimatortemperature025} for the estimator of the intensity considering this value of $h$. The criteria to minimize in order to find $\hat{h}$ is given in Figure \ref{criteriahyper} for the different values of $\alpha$. Similar results than \cite[Theorem3]{lacour16} can be derived in our context: the minimum penalty we can consider is achieved for $\alpha = 0$, and for $\alpha < 0$ the criteria leads to a value of $h$ closed to $h_{\min}$ with high probability and then induces a high bias. $\alpha = 0.25$ is positive but more likely to produce an under-biased estimator, which has happened in our case. The question of the choice of $\alpha$ remains open but a value of $\alpha$ since the asymptotical one, that is $\alpha = 1$, seems satisfying. This choice is also supported by \cite{varet17} (corresponding paper in progress) in the case of density estimation where numerical experiments are performed: the tuning parameter can be chosen equals to 1 without impact on the density estimation.  
 
\section{Numerical results}
\label{numericalresults}
In order to evaluate the performance of our estimation procedure, we present some simulation results. To be consistent with data, let us consider a model reproducing the temperature data and the spike times.

\medskip
 As in \cite{benth11}, we model the temperature $\theta_t$ as the sum of a trend seasonality function 
\[\Gamma_t = a + bt + c_1 \sin\left( \frac{2\pi t + \tau_1}{365 \times 24}\right) + c_2\sin\left( \frac{2\pi  t+\tau_2}{24}\right)\] 
corresponding to yearly and daily seasonality and a diffusion $X_t$ having dynamics
\[dX_t = -\vartheta X_tdt + \sigma dW_t\]
where $W_t$ is a standard Brownian motion. In \cite{benth11}, the temperature is modeled by a CARMA process with stochastic seasonal volatility but for simplicity we consider the simplest one corresponding to an Ornstein Uhlenbeck process. Using classical estimation procedures, we find $a = 12.06$, $b =  0.0000072$, $c_1 = 7.81$, $c_2 = -3.18$, $\tau_1 = -16924.50$, $\tau_2 = 10.84$, $\vartheta = 0.011 $, $\sigma = 0.46$.

\medskip
The spike intensity $\lambda_t$ is considered as an exponential function of the temperature
\[\lambda_t = a_0 \exp\left(a_1 \theta_t\right)\]
with $a_0 = 1033.8$ and $a_1 = -0.2$. As the quality of the estimation procedure depends on the interval of estimation, we consider three intervals: $\left[-1,29\right]$, $\left[-3,31\right]$ and $\left[-5,33\right]$. As for the estimation on data, we consider a local polynomial estimator of degree 1 with a Epanechnikov kernel. On each interval, we apply our estimation procedure and our bandwidth selection method with $h_{\min} = \frac{|I| \|K\|_1 \|K\|_{\infty}}{200}$ and $\mathcal{H} = \{h = h_{\min} + 0.1 i,\; h \leq 11 \}$ and we focus the analysis on the case where the tuning parameter $\alpha$ is set to 1. To evaluate the performance on our estimator, we consider the error 
\[e = \mathbb{E}\left(\frac{\int_I \left(q\left(x\right) - \hat{q}_{\hat{h}}\left(x\right)\right)^2dx}{\int_I q^2\left(x\right)dx}\right)\]
where $\hat{h}$ is the optimal bandwidth given by our estimation procedure and we compare it to the oracle error 
\[e_{o} = \underset{h \in \mathcal{H}}{\min}\; \mathbb{E}\left(\frac{\int_I \left(q\left(x\right) - \hat{q}_{h}\left(x\right)\right)^2dx}{\int_I q^2\left(x\right)dx} \right).\]
In practice, we consider estimators of these errors, $\hat{e}$ and $\hat{e}_o$. In a second time, we test if the intensity is of the form $a_0 \exp\left(a_1 x\right)$ and if it is constant, that is independent of the temperature.

\begin{figure}[h!]
             \centering
    \begin{subfigure}[b]{0.3\textwidth}
        \centering
        \includegraphics[width=\textwidth]{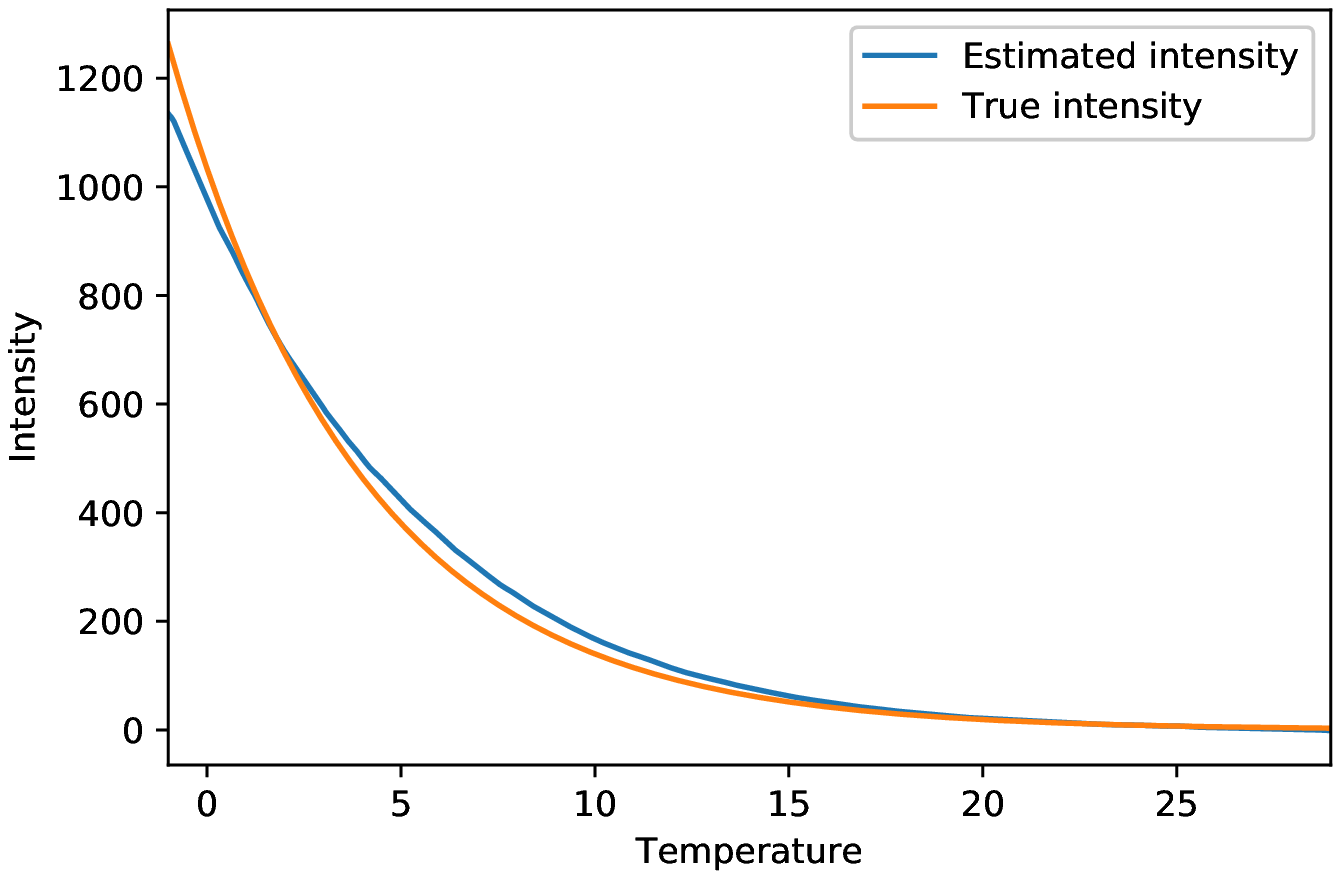}
        \caption{\it $I=\left[-1,29\right]$.}
    \end{subfigure}
    \begin{subfigure}[b]{0.3\textwidth}
        \centering
        \includegraphics[width=\textwidth]{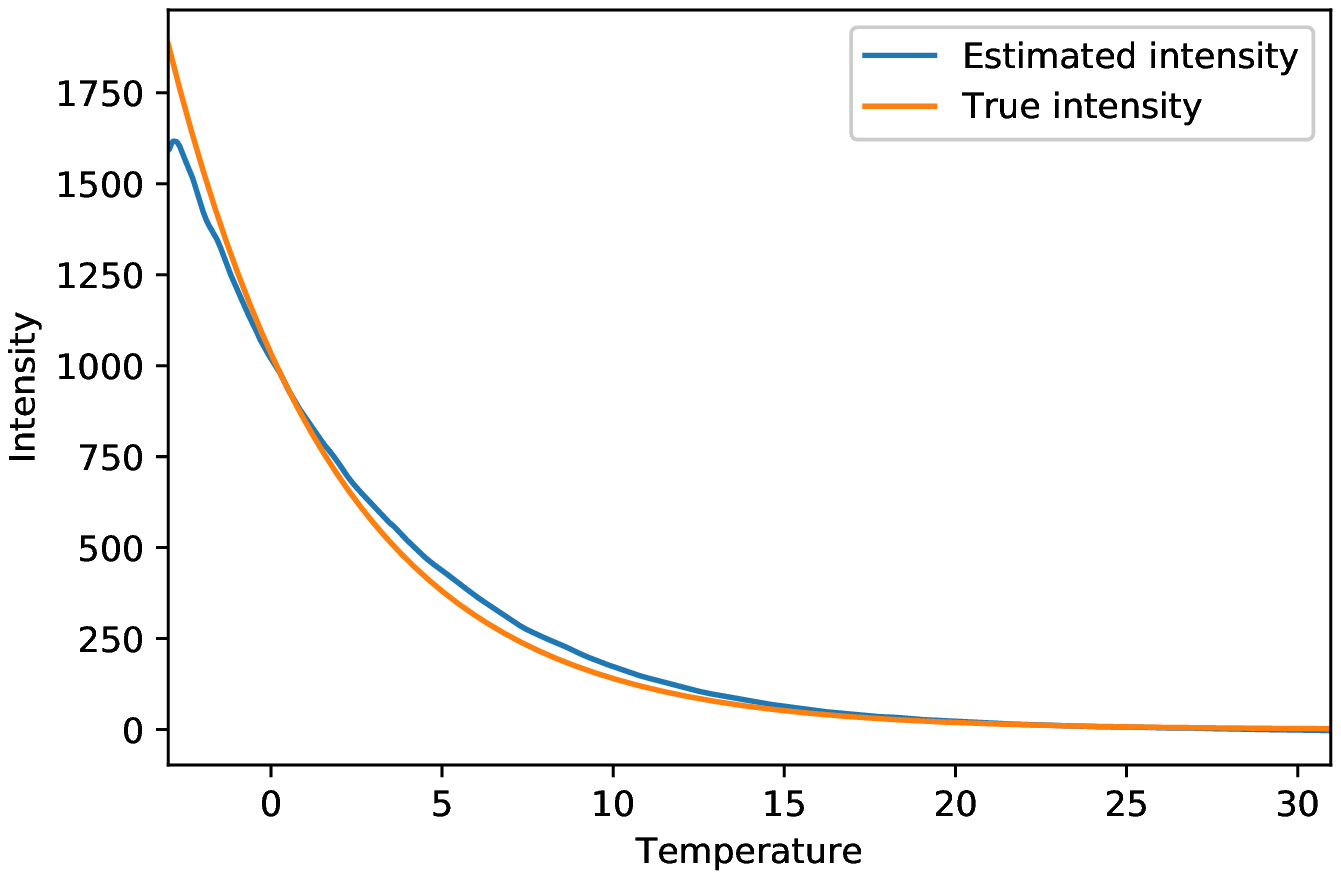}
        \caption{\it $I = \left[-3,31\right]$. }
    \end{subfigure}         
        \begin{subfigure}[b]{0.3\textwidth}
        \centering
        \includegraphics[width=\textwidth]{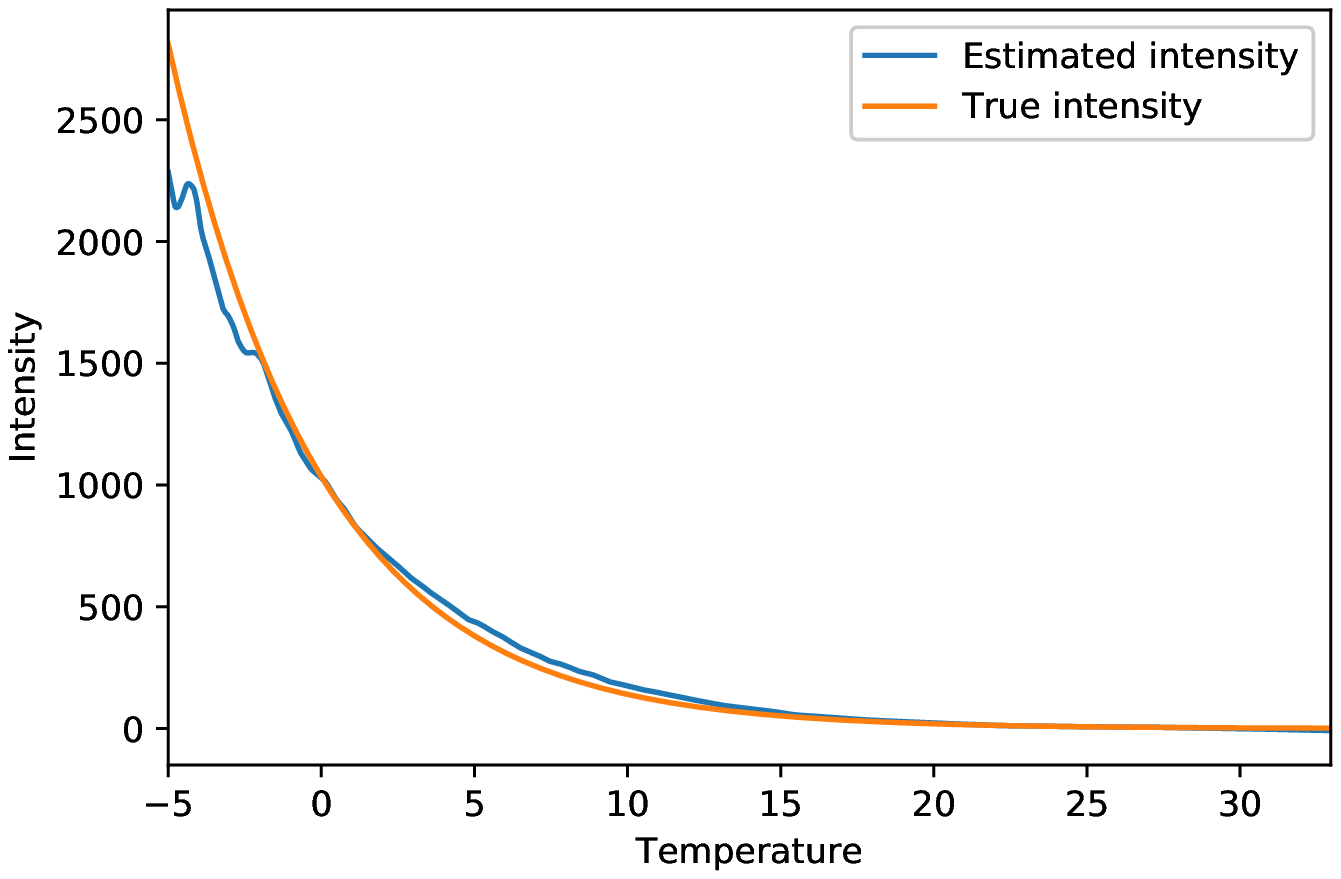}
        \caption{\it $I = \left[-5,33\right]$.}
    \end{subfigure}   
        \caption{\label{paramsimu} \it Mean of local polynomial estimators for different intervals with simulated data.}
    \end{figure}

\medskip
Results are given in Table \ref{simuresults} where 500 simulations of the model during $6$ years with a step time of one hour are considered. $\%$ converged corresponds to the percentage estimators that have been computed, meaning that their local time was large enough and the matrix $B$ invertible. This percentage diminishes with the length of the interval, and is very low for the last interval. However, this interval corresponds to the one we have estimated the parameters. One explanation is that the model does not capture all the features of the data. For instance, seasonal volatility is not modeled whereas it impacts the number of high and low values taken by the temperature. The two errors $\hat{e}$ and $\hat{e}_0$ increases with the length of the interval, caused by boundary effects: less values of the temperature are observed near the bounds. Furthermore, the ratio between $\hat{e}$ and $\hat{e}_0$ increases: the bandwidth selection procedure is less efficient for larger interval. This corresponds to the term $\frac{1}{\nu^2}$ in the oracle inequality. Columns $\%$ exponential and constant corresponds to the percentage of simulation for which the corresponding test has not been rejected at level $95\%$. Results are satisfying both for exponential and constant test. Estimators of $a_0$ and $a_1$ are consistent with the true parameters but present a small bias, probably due to the form of $M_n\left(\theta\right)$ that adds a weight term inside the norm. The mean of $\hat{q}_{\hat{h}}$ is represented in Figure \ref{paramsimu} for each interval $I$. One can see that there is a bias in the lower boundary for each interval. 

\begin{table}[h!]
\centering
\begin{scriptsize}
\begin{tabular}{|c|c|c|c|c|c|c|c|}
\hline
Interval  & $\hat{e}$ &  $\hat{e}_o$ & $\%$  converged & $\%$ exponential& $\%$ constant & $a_0$ & $a_1$   \\
\hline
$\left[-1,29\right]$ & 0.055 & 0.026 & 100 & 97 & 0 &$ \left[1030.60,1059,65\right] $& $\left[-0.210,-0.202\right] $	\\
\hline
$\left[-3,31\right]$ & 0.082 & 0.03 & 81 & 97 & 0.25 & $\left[1015.67,1041.91\right]$& $\left[-0.223,-0.202\right]$ 	\\
\hline
$\left[-5,33\right]$ &0.19 & 0.04 & 12.6 & 73 & 0 & $\left[935.77,973.50\right]$ & $\left[-0.253,-0.2212\right]$ \\
\hline
\end{tabular}
 \end{scriptsize}
 \caption{\label{simuresults} \it Performance of the local polynomial estimation procedure and parametrical test on different intervals for simulated data.} 
 \end{table}

\section{Proof of Proposition \ref{oracle}}
\label{prooforacle}
In order to prove Proposition \ref{oracle}, we need Proposition \ref{biasvariance} and Proposition \ref{oraclewithpen}. Proposition \ref{biasvariance} gives an approximation of the error by the bias and the variance. This proposition is similar to the one of \cite[Proposition 4.1]{lerasle16} in the context of density estimation.

\medskip
During the proof, $\tilde{C}$ denotes a constant that can change from line to line. $\tilde{C}\left(\cdot\right)$ denotes a constant depending on $\cdot$ that can also change from line to line.

\begin{proposition} \label{biasvariance} Assume \ref{assumptionsX} and \ref{assumptionkernel}. Let $x \geq 1$, $\eta \in \left(0,1\right]$. Let $\mathcal{H}$ a finite subset of $\left(0,\infty\right)$ such that $\min \mathcal{H} = h_{\min} \geq \frac{\|K\|_{\infty} \|K\|_1 |I|}{n}$ and $\max \mathcal{H} \leq  \frac{2}{3} \frac{| I |}{\Delta} $. Let $\hat{q}_h$ the local polynomial estimator defined in \eqref{estimator} and $\hat{h}$ the bandwidth defined in \eqref{optimalh}. With conditional probability given $\mathcal{F}^X_T$ larger than $1- \tilde{C} |\mathcal{H}|e^{-x}$, on the event $D\left(I,\nu\right)$, for any $h \in \mathcal{H}$,
\[\|q-\hat{q}_{h}\|_I^2  \leq \left(1+\eta\right)\left(\|q- q_{h}\|_I^2 + V_h \right) + \tilde{C} \frac{|I|\left(1 + \|q\|_{I,\infty} \|l_T\|_{I,\infty} |I| \| K\|_1^2\right)  x^2}{n A^2_K \nu^2 T^2 \eta^3} \]
and
\[ \|q- q_{h}\|_I^2 + V_{h}  \leq \left(1+\eta\right) \|q- \hat{q}_{h}\|_I^2 + \tilde{C} \frac{|I|\left(1 + \|q\|_{I,\infty}\|l_T\|_{I,\infty}\| |I| K\|_1^2\right) x^2 }{n A^2_K \nu^2 T^2 \eta^3} \]
where 
\[V_h =\frac{1}{n} \int_{0}^T \left(\int_{I} \left(w\left(x,h,\frac{X_s-x}{h}\right) K_h\left(X_s-x\right)\right)^2 {\bf 1}_{X_s \in I} dx\right) q\left(X_s\right)ds\]
and $A_K$ is a constant depending on $K$ which is introduced in Proposition \ref{B}.
\end{proposition}

\begin{proposition} \label{oraclewithpen}Assume \ref{assumptionsX} and \ref{assumptionkernel}. Let $x \geq 1$, $\theta \in \left(0,1\right)$. Let $\mathcal{H}$ a finite subset of $\left(0,\infty\right)$ such that $\min \mathcal{H} = h_{\min} \geq \frac{\|K\|_{\infty} \|K\|_1 |I|}{n}$ and $\max \mathcal{H} \leq  \frac{2}{3} \frac{| I |}{\Delta} $. Let $\hat{q}_h$ the local polynomial estimator defined in \eqref{estimator} and $\hat{h}$ the bandwidth defined in \eqref{optimalh}. Let $\hat{V}_{h,h_{\min}}$ defined in Equation \eqref{Vhhmin} and $pen\left(\alpha\right)$ defined in Equation \eqref{pen} with $\alpha > 0$. With conditional probability given $\mathcal{F}^X_T$ larger than $1-C_1|\mathcal{H}|e^{-x}$, on the event $D\left(I,\nu\right)$, for any $h \in \mathcal{H}$,
\begin{align*}\left(1-\theta\right)\|\hat{q}_{\hat{h}} - q\|_I^2 &\leq \left(1+\theta\right) \| \hat{q}_h - q\|_I^2 + \left(\text{pen}_{\alpha}\left(h\right) - 2\hat{V}_{h,h_{\min}}\right) -  \left(\text{pen}_{\alpha}\left(\hat{h}\right) - 2\hat{V}_{\hat{h},h_{\min}}\right) \\
&+ \frac{C_2}{\theta}\|q_{h_{\min}}-q\|_I^2 + \frac{C\left(K\right) |I|}{\nu^2 T^2 \theta}\left( \frac{\left(1+\|q\|_{I,\infty}\|l_T\|_{I,\infty} |I|\right) x^2}{n} + \frac{x^3 |I|}{n^2 h_{\min}}\right) 
\end{align*}
where $C_1$ and $C_2$ are constant and $C\left(K\right)$ is a constant depending on $K$.
\end{proposition}

\subsection{Proof of Proposition \ref{biasvariance}}

In the following, we work conditionally on $\mathcal{F}_T^X$. We also work on the event $D\left(I,\nu\right)$. Conditionally on $\mathcal{F}_T^X$, the process $N$ is a inhomogeneous Poisson process. Proposition \ref{bernstein} and Proposition \ref{ustatistic} are then verified taking the conditional probability given $\mathcal{F}_T^X$.

\medskip
The norm $\|q-\hat{q}_h\|_I^2$ is the sum of $\|q-q_{h}\|_I^2$ which is a bias term, $\|\hat{q}_{h} - q_{h}\|_I$ which is a variance term and the cross term $2<q-q_{h},q_{h}-\hat{q}_{h}>_I$. In order to control $\|q-\hat{q}_{h}\|_I^2$ by $\|q-q_{h}\|_I^2 + V_h$, we will control the variance term by $V_h$ and the cross term by $\|q-q_{h}\|_I^2 + V_{h}$. We consider a real number $x \geq 1$ in the following.

\medskip
\paragraph{{\bf Control of the variance term}}  First, let us control the term $\|q_{h} - \hat{q}_{h}\|_I^2$. This term is equal to 
\[\frac{1}{n^2}\int_{I} \left(\int_{0}^T w\left(x,h,\frac{X_s-x}{h}\right) K_{h}\left(X_s-x\right){\bf 1}_{X_s \in I}dM_s\right)^2 dx\] which can be written as the sum of 
\begin{equation} \label{decompositionvariance1}
\frac{1}{n^2}\int_{0}^T \int_{I} w^2\left(x,h,\frac{X_s-x}{h}\right) K^2_h\left(X_s-x\right){\bf 1}_{X_s \in I} dx dN_s
\end{equation}
and 
\begin{equation} \label{decompositionvariance2}
\frac{2}{n^2} \hspace{-0.5em} \int_{0}^T \hspace{-0.5em} \int_{0}^{s^-} \hspace{-0.9em}\int_{I} w\left(x,h,\frac{X_s-x}{h}\right)w\left(x,h,\frac{X_u-x}{h}\right)K_h\left(X_u-x\right)K_h\left(X_s-x\right) dx {\bf 1}_{X_s \in I}{\bf 1}_{X_u \in I} dM_u dM_s. 
\end{equation}
The term \eqref{decompositionvariance1} is a simple Poisson integral and can be controlled with Proposition \ref{bernstein}. 
As 
\begin{align*}
\frac{1}{n^2}\int_{I} w^2\left(x,h,\frac{X_s-x}{h}\right)K^2_h\left(X_s-x\right){\bf 1}_{X_s \in I}dx &\leq \frac{\|K\|_1 \|K\|_{\infty} |I|^2}{n^2 A^2_K \nu^2 T^2 h}\\
&\leq \frac{|I|}{n A^2_K \nu^2 T^2}
\end{align*}
and 
\[\frac{1}{n^4} \int_{0}^T \left(\int_{I} w^2\left(x,h,\frac{X_s-x}{h}\right) K^2_h\left(X_s-x\right){\bf 1}_{X_s \in I}dx\right)^2 d\Lambda_s \]
is bounded by 
\begin{align*}
\frac{1}{n^2}\underset{z \in I}{\sup} \int_{I} w^2\left(x,h,\frac{z-x}{h}\right) K^2_h\left(z-x\right) dx V_h&\leq \frac{ \|K\|_1 \|K\|_{\infty}V_h |I|^2}{A^2_K \nu^2 T^2 h n^2} \\
&\leq \frac{ V_h |I|}{n A^2_K \nu^2 T^2},
\end{align*}
with conditional probability given $\mathcal{F}_T^X$ larger than $1-2|\mathcal{H}| e^{-x}$, on $D\left(I,\nu\right)$, for any $h \in \mathcal{H}$,
\[ | \frac{1}{n^2}\int_{0}^T \int_{I} w^2\left(x,h,\frac{X_s-x}{h}\right)K^2_h\left(X_s-x\right) dx {\bf 1}_{X_s \in I} dN_s - V_h | \leq \sqrt{ \frac{2|I|x}{n A^2_K \nu^2 T^2} V_h} +  \frac{|I| x}{3 n A^2_K \nu^2 T^2}.\]
Using Young's inequality 
\begin{equation} \label{younginequality}2ab \leq \epsilon a^2+ \frac{b^2}{\epsilon}, \text{ for all } \epsilon > 0,
\end{equation}
we find, for $\theta > 0$, with conditional probability given $\mathcal{F}_T^X$ larger than $1-2|\mathcal{H}| e^{-x}$, on $D\left(I,\nu\right)$, for any $h \in \mathcal{H}$,
\begin{equation} \label{variancepart1}
| \frac{1}{n^2} \int_{0}^T \int_{I}w^2\left(x,h,\frac{X_s-x}{h}\right) K^2_h\left(X_s-x\right) dx {\bf 1}_{X_s \in I}dN_s - V_h | \leq \theta V_h +  \tilde{C} \frac{x |I|}{n A^2_K \nu^2 T^2 \theta}.
\end{equation}

The term \eqref{decompositionvariance2} is an U-statistics which can be controlled with Proposition \ref{ustatistic}. With conditional probability given $\mathcal{F}^X_T$  larger than $1-6.44 e^{-x}$, it is dominated in absolute value and on $D\left(I,\nu\right)$, by 
\[\tilde{C} \left(C\sqrt{x} + Dx + Bx^{\frac{3}{2}} + Ax^{2}\right)\]
with $A$, $B$, $C$ and $D$ defined in the following.
We have 
\begin{align*}A &= \frac{1}{n^2}\underset{\left(u,s\right) \in \left[0,T\right]^2}{\sup} \int_{I} w\left(x,h,\frac{X_s-x}{h}\right)w\left(x,h,\frac{X_u-x}{h}\right) K_h\left(X_u-x\right)K_h\left(X_s-x\right){\bf 1}_{X_s \in I}{\bf 1}_{X_u \in I}dx \\
&\leq \frac{\|K\|_{\infty} \|K\|_1 |I|^2}{n^2 A^2_K h\nu^2 T^2}\\
&\leq \frac{|I|}{n A^2_K \nu^2 T^2}.
\end{align*} 

The term 
\[
\begin{split}
B^2 &=  \max\{\underset{s \leq T}{\sup} \int_{0}^s \left(\int_{I} w\left(x,h,\frac{X_s-x}{h}\right)w\left(x,h,\frac{X_u-x}{h}\right) K_h\left(X_u-x\right)K_h\left(X_s-x\right){\bf 1}_{X_s \in I}{\bf 1}_{X_u \in I}dx\right)^2 d\Lambda_u,\\
& \underset{u \leq T}{\sup} \int_{u}^T \left(\int_{I} w\left(x,h,\frac{X_s-x}{h}\right)w\left(x,h,\frac{X_u-x}{h}\right) K_h\left(X_u-x\right)K_h\left(X_s-x\right){\bf 1}_{X_s \in I}{\bf 1}_{X_u \in I}dx\right)^2 d\Lambda_s \}
\end{split}
\]
is bounded by, using the occupation time formula, 
\begin{align*}
&\frac{ |I|^4 \underset{\left(v,z\right)\in I^2}{\sup} \int_{I} K_h\left(z-x\right)K_h\left(v-x\right)dx  \underset{u \in \left[0,T\right]}{\sup} \int_I \int_{I} K_h\left(z-x\right) K_h\left(X_u-x\right){\bf 1}_{X_u \in I} dx l_T^z n q\left(z\right) dz}{n^4 A^4_K \nu^4 T^4} \\
&\leq \frac{\|K\|_{\infty} \|K\|^3_1 \|q\|_{I,\infty} \| l_T\|_{I,\infty}|I|^4}{h n^3 A_K^4 \nu^4 T^4}\\
&\leq \frac{\|K\|^2_1 \|q\|_{I,\infty} \| l_T\|_{I,\infty} |I|^3}{n^2 A^4_K \nu^4 T^4}.
\end{align*}
Using again the occupation time formula and Young's inequality for convolutions,
\[C^2 = \int_{0}^T \int_{0}^s \left(\int_{I} w\left(x,h,\frac{X_s-x}{h}\right)w\left(x,h,\frac{X_u-x}{h}\right) K_h\left(X_u-x\right)K_h\left(X_s-x\right){\bf 1}_{X_s \in I}{\bf 1}_{X_u \in I}dx\right)^2 d\Lambda_u d\Lambda_s\] is bounded by
\begin{align*}
& \frac{ \|l_T\|_{I,\infty} \| q  \|_{I,\infty} |I|^2}{n^3 A^2_K \nu^2 T^2} \int_0^T \int_{I} \left(\int_{I} w\left(x,h,\frac{X_s-x}{h}\right) K_h\left(X_s-x\right){\bf 1}_{X_s \in I} K_h\left(z-x\right) dx\right)^2  dz d\Lambda_s\\
&\leq \frac{ \|l_T\|_{I,\infty} \| q  \|_{I,\infty} |I|^2}{n^3 A^2_K \nu^2 T^2} \int_{0}^T \| w\left(\cdot,h,\frac{X_s-\cdot}{h}\right)K_h\left(X_s-\cdot\right){\bf 1}_{X_s \in I} {\bf 1}_{\cdot \in I} \ast K_h\left( \cdot \right) \|^2 d\Lambda_s\\
&\leq  \frac{ \|l_T\|_{I,\infty} \| q  \|_{I,\infty} |I|^2}{n^3 A^2_K \nu^2 T^2} \int_{0}^T \| w\left(\cdot,h,\frac{X_s-\cdot}{h}\right) K_h\left(X_s-\cdot\right){\bf 1}_{X_s \in I} {\bf 1}_{\cdot \in I}    \|^2 \| K_h\|^2_1 d\Lambda_s \\
&=  \frac{\|l_T\|_{I,\infty} \| q  \|_{I,\infty} \|K\|_1^2 V_h |I|^2}{n A^2_K \nu^2 T^2}.
\end{align*}
Using twice Cauchy Schwarz inequality, $D$ which is equal to 
\[\mathop{\underset{{\int_{0}^T a^2_u d\Lambda_u = 1,}}{\sup}}_{\int_{0}^T b^2_s d\Lambda_s = 1}\; \int_{0}^T a_u \int_{u}^T b_s \int_{I} w\left(x,h,\frac{X_s-x}{h}\right)w\left(x,h,\frac{X_u-x}{h}\right) K_h\left(X_u-x\right)K_h\left(X_s-x\right){\bf 1}_{X_s \in I}{\bf 1}_{X_u \in I}dx d\Lambda_s d\Lambda_u \]
can be bounded in the same way than $C$ and  
\[D \leq \sqrt{ \frac{\|l_T\|_{I,\infty} \| q  \|_{I,\infty} \|K\|_1^2 V_h |I|^2}{n  A^2_K \nu^2 T^2}}.\]
After applying \eqref{younginequality}, we obtain for $\theta \in \left(0,1\right]$ that with conditional probability given $\mathcal{F}_T^X$ larger than $1-\tilde{C} |\mathcal{H}|e^{-x}$, on $D\left(I,\nu\right)$,  for any $h \in \mathcal{H}$,
\begin{align*} \label{variancepart2}
&\frac{2}{n^2}|\int_{0}^T\hspace{-0.5em} \int_{0}^{s^-}\hspace{-0.5em} \int_{I}w\left(x,h,\frac{X_u-x}{h}\right)w\left(x,h,\frac{X_s-x}{h}\right) K_h\left(X_u-x\right){\bf 1}_{X_u \in I}K_h\left(X_s-x\right){\bf 1}_{X_s \in I}dx dM_u dM_s|  \\ \numberthis
&\leq  \theta V_h  + \tilde{C} \frac{|I|\left(1+\|q\|_{I,\infty} \|l_T\|_{I,\infty} |I| \|K\|_1^2\right)x^2}{\nu^2 T^2 A^2_K n \theta}. 
\end{align*}
Combining \eqref{variancepart1} and \eqref{variancepart2} with $\tilde{\theta} = \frac{\theta}{2}$, we find for $\tilde{\theta} \in \left(0,1\right]$, with conditional probability given $\mathcal{F}_T^X$ larger than $1-\tilde{C} |\mathcal{H}|e^{-x}$, on $D\left(I,\nu\right)$, for any $h \in \mathcal{H}$,
\begin{equation} \label{controlvariance}
| \|q_{h} - \hat{q}_{h}\|_I^2 - V_h | \leq \tilde{\theta} V_h + \tilde{C} \frac{|I|\left(1+\|q\|_{I,\infty}\|l_T\|_{I,\infty}|I| \|K\|_1^2\right)x^2}{n A^2_K \nu^2 T^2 \tilde{\theta}}.
 \end{equation}
%
%
\medskip
\paragraph{{\bf Control of the cross term}} The cross term $<\hat{q}_	{h} - q_{h}, q_{h} - q>_I$ is equal to 
\[\frac{1}{n}\int_{0}^T \left(\int_{I} \left(q_{h}\left(x\right) - q\left(x\right)\right) w\left(x,h,\frac{X_s-x}{h}\right) K_h\left(X_s-x\right){\bf 1}_{X_s \in I} dx\right) dM_s\]
and thus can be controlled using Proposition \ref{bernstein}. Using Cauchy Schwarz,
\begin{align*}
\frac{1}{n}\int_{I} \left(q_{h}\left(x\right) - q\left(x\right)\right) w\left(x,h,\frac{X_s-x}{h}\right) K_h\left(X_s-x\right){\bf 1}_{X_s \in I}dx &\leq \frac{1}{n}\| q_{h} - q \|_I \sqrt{ \frac{ \|K\|_1 \|K\|_{\infty} |I|^2}{h A^2_K \nu^2 T^2}}\\    
&\leq  \| q_{h} - q \|_I \sqrt{\frac{|I|}{n A^2_K \nu^2 T^2 } }
\end{align*}
and for $\theta > 0$, $\frac{1}{n}\int_{I} \left(q_{h}\left(x\right) - q\left(x\right)\right) w\left(x,h,\frac{X_s-x}{h}\right) K_h\left(X_s-x\right){\bf 1}_{X_s \in I}dx$ is bounded by
\begin{equation} \label{crosstermbound1}
\frac{\theta \|q_{h}-q\|_I^2}{2} + \frac{|I|x^2}{2n\theta A^2_K \nu^2 T^2}.
\end{equation}
Using Cauchy Schwarz inequality again, 
\[\frac{1}{n^2}\int_{0}^T \left( \int_{I} \left(q_{h}\left(x\right) - q\left(x\right)\right) w\left(x,h,\frac{X_s-x}{h}\right) K_h\left(X_s-x\right){\bf 1}_{X_s \in I} dx   \right)^2 d\Lambda_s \leq \|q_{h} - q\|_I^2 V_h.\]
We can also bound this term using the occupation time formula and Young's inequality for convolution by 
\[ \frac{|I|^2\|l_T\|_{I,\infty} \|q\|_{I,\infty}}{n A^2_K \nu^2 T^2} \| K_h\left(\cdot\right) \ast \left(q_{h}-q\right)\left( \cdot\right) {\bf 1}_{\cdot \in I} \|^2 \leq \frac{ \|l_T\|_{I,\infty} \|q\|_{I,\infty} \| K\|_1^2 |I|^2}{n A^2_K \nu^2 T^2} \|q_{h}-q\|_I^2\]
Thus, 
\[\frac{1}{n^2}\int_{0}^T \left( \int_{I} \left(q_{h}\left(x\right) - q\left(x\right)\right) w\left(x,h,\frac{X_s-x}{h}\right)  K_h\left(X_s-x\right){\bf 1}_{X_s \in I} dx   \right)^2 d\Lambda_s\]
is bounded by 
\[\|q-q_{h}\|_I^2 \sqrt{V_{h}} \sqrt{\frac{ \|l_T\|_{I,\infty} \|q\|_{I,\infty} \| K\|_1^2 |I|^2 }{n A^2_K \nu^2 T^2}}\]
and for $\theta, u > 0$,
\[\sqrt{\frac{2}{n^2} \int_{0}^T \left( \int_{I} \left(q_{h}\left(x\right) - q\left(x\right)\right)w\left(x,h,\frac{X_s-x}{h}\right)  K_h\left(X_s-x\right){\bf 1}_{X_s \in I} dx   \right)^2 d\Lambda_s x}\]
is bounded by
\[ 
\theta \|q-q_{h}\|_I^2 + \frac{1}{2\theta} \sqrt{V_h} \sqrt{\frac{ \|l_T\|_{I,\infty} \|q\|_{I,\infty}  \| K\|_1^2 |I|^2 x}{n A^2_K n \nu^2 T^2}} \leq \theta \|q-q_{h}\|_I^2 + \frac{uV_h}{\theta} + \frac{x \|l_T\|_{I,\infty}   \|q\|_{I,\infty}  \| K\|_1^2 |I|^2}{16 \theta u A^2_K n  \nu^2 T^2}.\]
Thus, for $\theta \in \left(0,1\right]$, with $u = \theta^2$, we obtain, with conditional probability given $\mathcal{F}^X_T$ larger than $1-\tilde{C} |\mathcal{H}|e^{-x}$, on $D\left(I,\nu\right)$, for any $h \in \mathcal{H}$,
\[\sqrt{\frac{2}{n^2} \int_{0}^T \left( \int_{I} \left(q_{h}\left(x\right) - q\left(x\right)\right) w\left(x,h,\frac{X_s-x}{h}\right)K_h\left(X_s-x\right){\bf 1}_{X_s \in I}dx   \right)^2 d\Lambda_s x}\]
is bounded by 
\begin{equation} \label{crosstermbound2}
\theta\left(\|q-q_{h}\|_I^2 + V_h\right) + \tilde{C} \frac{x \|l_T\|_{I,\infty} \|q\|_{I,\infty} \| K\|_1^2 |I|^2 }{\theta^3 n A^2_K \nu^2 T^2}.
\end{equation}
Using Proposition \ref{bernstein}, bounds \eqref{crosstermbound1} and \eqref{crosstermbound2}, we have for $\theta \in\left(0,1\right]$, with conditional probability given $\mathcal{F}^X_T$ larger than $1-\tilde{C} |\mathcal{H}|e^{-x}$, on $D\left(I,\nu\right)$, for any $h \in \mathcal{H}$,
\begin{equation} \label{controlcross}
2|<q_{h}-\hat{q}_{h},q_{h} -q>_I| \leq 3\theta \left(\|q-q_{h}\|_I^2 + V_h\right) + \tilde{C}  \frac{|I| \left(1+\|q\|_{I,\infty}\|l_T\|_{I,\infty}|I| \| K\|_1^2\right) x^2}{n A^2_K\nu^2 T^2 \theta^3}.\end{equation}
Combining \eqref{controlvariance} and \eqref{controlcross}, we find for $\theta \in\left(0,1\right]$, with conditional probability given $\mathcal{F}^X_T$ larger than $1-\tilde{C} |\mathcal{H}|e^{-x}$, on $D\left(I,\nu\right)$, for any $h \in \mathcal{H}$,
\begin{equation} \label{endresult1}
\|q-\hat{q}_{h}\|^2_I - \|q-q_h\|^2_I \leq 3\theta \|q-q_{h}\|_I^2  + \left(1+4\theta\right)V_h + \tilde{C}  \frac{|I| \left(1+\|q\|_{I,\infty}\|l_T\|_{I,\infty} |I| \| K\|_1^2\right) x^2}{n A^2_K\nu^2 T^2 \theta^3}
\end{equation}
and 
\begin{equation} \label{endesult2}\|q-\hat{q}_{h}\|^2_I - \|q-q_h\|^2_I \geq -3\theta \left(\|q-q_{h}\|_I^2 + V_h\right) + \left(1-\theta\right)V_h -  \tilde{C}  \frac{|I|\left(1+\|q\|_{I,\infty}\|l_T\|_{I,\infty}|I| \| K\|_1^2\right) x^2}{n A^2_K\nu^2 T^2 \theta^3}.
\end{equation}
Taking $\theta = \frac{\eta}{4}$ for \eqref{endresult1} and $\theta = \frac{\eta}{4\left(1+\eta\right)}$ for \eqref{endresult2} with $\eta \in \left(0,1\right]$, the proof is achieved.

\subsection{Proof of Proposition \ref{oraclewithpen}}

We continue to work conditionally on $\mathcal{F}^X_T$ and on the event $D\left(I,\nu\right)$.
The beginning of the proof is similar to the one of \cite[Theorem 9]{lacour16}. For any $h \in \mathcal{H}$, 
\begin{align*}
\| \hat{q}_{\hat{h}}-q\|_I^2 + \text{pen}_{\alpha}\left(\hat{h}\right) &= \|\hat{q}_{\hat{h}}-\hat{q}_{h_{\min}}\|_I^2 + \text{pen}_{\alpha}\left(\hat{h}\right) + \|\hat{q}_{h_{\min}}-q\|_I^2 + 2<\hat{q}_{\hat{h}}-\hat{q}_{h_{\min}},\hat{q}_{h_{\min}}-q>_I\\
&\leq \|\hat{q}_{h} - \hat{q}_{h_{\min}}\|_I^2 + \text{pen}_{\alpha}\left(h\right) + \|\hat{q}_{h_{\min}}-q\|_I^2 + 2<\hat{q}_{\hat{h}}-\hat{q}_{h_{\min}},\hat{q}_{h_{\min}}-q>_I\\
&\leq \|\hat{q}_{h}-q\|_I^2 + 2\|q-\hat{q}_{h_{\min}}\|_I^2 + 2<\hat{q}_{h}-q,q-\hat{q}_{h_{\min}}>_I + \text{pen}_{\alpha}\left(h\right)\\
&+2<\hat{q}_{\hat{h}}-\hat{q}_{h_{\min}},\hat{q}_{h_{\min}}-q>_I.
\end{align*}
We then have 
\begin{equation}\label{controlpen}
\begin{split}
\| \hat{q}_{\hat{h}}-q\|_I^2 \leq \| \hat{q}_{h}-q\|_I^2 &+ \left( \text{pen}_{\alpha}\left(h\right) - 2<\hat{q}_{h}-q,\hat{q}_{h_{\min}}-q>_I\right)\\
&-\left( \text{pen}_{\alpha}\left(\hat{h}\right) - 2<\hat{q}_{\hat{h}}-q,\hat{q}_{h_{\min}}-q>_I\right).  
\end{split}
\end{equation}
We want to approach $<\hat{q}_{h}-q,\hat{q}_{h_{\min}}-q>_I = S\left(h,h_{\min} \right)$ by $\hat{V}_{h,h_{\min}}$ with $S\left(h,h' \right) = <\hat{q}_{h}-q_{h},q_{h'}-q>_I$ for $h,\; h' \in \mathcal{H}$.
We have 
\begin{align*}
<\hat{q}_{h}-q,\hat{q}_{h_{\min}}-q>_I &= <\hat{q}_{h}-q_{h}+q_{h}-q,\hat{q}_{h_{\min}}-q_{h_{\min}}+q_{h_{\min}}-q>_I\\
&= <\hat{q}_{h}-q_{h},\hat{q}_{h_{\min}}-q_{h_{\min}}>_I + S\left(h,h_{\min}\right) + S\left(h_{\min},h\right) \\
&+ <q_{h}-q,q_{h_{\min}}-q>_I.
\end{align*}
Furthermore, we can easily show that 
\[
<\hat{q}_{h}-q_{h},\hat{q}_{h_{\min}}-q_{h_{\min}}>_I = \hat{V}_{h,h_{\min}} + U\left(h,h_{\min}\right)
\]
with \begin{align*}
\hat{V}_{h,h_{\min}}  = \\
&\frac{1}{n^2} \int_{0}^T \int_{I} w\left(x,h,\frac{X_s-x}{h}\right) w\left(x,h_{\min},\frac{X_s-x}{h_{\min}}\right)  K_h\left(X_s-x\right)K_{h_{\min}}\left(X_s-x\right){\bf 1}_{X_s \in I}dx dN_s,
\end{align*}
\[U\left(h,h_{\min}\right) = \frac{1}{n^2}\int_0^T \int_0^{s^{-}} \left(G_{h,h_{\min}}\left(X_u,X_s\right)+G_{h_{\min},h}\left(X_s,X_u\right)\right)dM_u dM_s\]
and 
\[\begin{split}
&G_{h,h'}\left(X_u,X_s\right) =\\
& \int_{I} w\left(x,h,\frac{X_u-x}{h}\right) w\left(x,h_{\min},\frac{X_s-x}{h_{\min}}\right) K_h\left(X_u-x\right){\bf 1}_{X_u \in I}K_{h_{\min}}\left(X_s-x\right){\bf 1}_{X_s \in I}dx.
\end{split}
\]
Thus, 
\begin{equation} \label{decompositionpen}
<\hat{q}_{h}-q,\hat{q}_{h_{\min}}-q>_I = \hat{V}_{h,h_{\min}}  + U\left(h,h_{\min}\right) + S\left(h,h_{\min}\right) + S\left(h_{\min},h\right) + <q_{h}-q,q_{h_{\min}}-q>_I.
\end{equation}

\medskip
In the following, we consider two real numbers $x \geq 1$ and $\theta' \in \left(0,1\right)$.

\medskip
\paragraph{{\bf Control of $U\left(h,h_{\min}\right)$}} The term $U\left(h,h_{\min}\right)$ can be controlled using Proposition \ref{ustatistic}. With conditional probability given $\mathcal{F}^X_T$ larger than $1-\tilde{C} |\mathcal{H}| e^{-x}$, on $D\left(I,\nu\right)$, for any $h \in \mathcal{H}$, 
\[ |U\left(h,h_{\min}\right)| \leq \tilde{C} \left(C \sqrt{x} + Dx + B x^{\frac{3}{2}} + Ax^2\right)\]
with $A$, $B$, $C$ and $D$ defined below.
We have 
\begin{align*}
A &= \frac{1}{n^2}\underset{\left(u,v\right) \in I^2}{\sup}\; |G_{h,h_{\min}}\left(u,v\right) + G_{h_{\min},h}\left(v,u\right) |\\
&\leq \frac{2}{n^2}\underset{\left(u,v\right) \in I^2}{\sup}\;  |G_{h,h_{\min}}\left(u,v\right)|\\
&\leq \frac{2 \|K\|_{\infty} \|K\|_1 |I|^2}{A^2_K \nu^2 T^2 n^2 h_{\min}}.
\end{align*}
We have $B^2$ equal to 
\[\frac{1}{n^4}\max\{\underset{s \leq T}{\sup} \int_{0}^s \left(G_{h,h_{\min}}\left(X_u,X_s\right)+G_{h,h_{\min}}\left(X_s,X_u\right)\right)^2 d\Lambda_u, \underset{u \leq T}{\sup} \int_{u}^T \left(G_{h,h_{\min}}\left(X_u,X_s\right)+G_{h,h_{\min}}\left(X_s,X_u\right)\right)^2 d\Lambda_s   \}\]
and bounded by
\[
\begin{split}
&\leq \frac{1}{n^4} \underset{u \leq T}{\sup} \int_0^T \left(G_{h,h_{\min}}\left(X_u,X_s\right)+G_{h,h_{\min}}\left(X_s,X_u\right)\right)^2 d\Lambda_s \\
&\leq \frac{2}{n^4} \underset{u \leq T}{\sup} \int_0^T \left(G_{h,h_{\min}}\left(X_u,X_s\right)\right)^2 d\Lambda_s+\int_0^T \left(G_{h,h_{\min}}\left(X_s,X_u\right)\right)^2 d\Lambda_s.
\end{split}
\]
Using Cauchy Schwarz inequality, 
\begin{equation} \label{UB1}
\frac{1}{n^4}\int_{0}^T \left(G_{h,h_{\min}}\left(X_s,X_u\right)\right)^2 d\Lambda_s \leq \frac{V_h \|K\|^2 |I|^2}{n^2 A^2_K \nu^2 T^2 h_{\min}}.
\end{equation}
Using the occupation time formula and Young's inequality for convolutions,
\begin{align*} 
\frac{1}{n^4}\int_{0}^T \left(G_{h,h_{\min}}\left(X_u,X_s\right) \right)^2 d\Lambda_s &\leq 
\frac{ \|q\|_{I,\infty} \| l_T \|_{I,\infty}|I|^4 \|K_h\left(X_u-\cdot\right) \ast K_{h_{\min}}\left(\cdot\right) \|^2}{n^3 A^4_K \nu^4 T^4}\\
&\leq \frac{ \|q\|_{I,\infty} \|l_T\|_{I,\infty} \|K\|_1^2  \| K \|^2 |I|^4}{h_{\min} n^3 A^4_K \nu^4 T^4}. \numberthis \label{UB2}
\end{align*}
Thus, combining \eqref{UB1} and \eqref{UB2},
\[B^2 \leq  \frac{2 V_h \|K\|^2 |I|^2}{A^2_K n^2 \nu^2 T^2 h_{\min}} + \frac{2 \|q\|_{I,\infty} \|l_{T}\|_{I,\infty} \|K\|_1^2  \| K \|^2|I|^4}{h_{\min} A^4_K n^3 \nu^4 T^4},\]
\[B \leq \sqrt{\frac{2V_h \|K\|^2 |I|^2}{A^2_K n^2 \nu^2 T^2 h_{\min}}} + \sqrt{\frac{ 2\|q\|_{I,\infty} \|l_{T}\|_{I,\infty}|I|^4 \|K\|_1^2  \| K \|^2}{n^3 A^4_K h_{\min} \nu^4 T^4}}\]
and 
\begin{align*}
Bx^{\frac{3}{2}} &\leq \frac{\theta'}{3} V_h + \frac{3}{2\theta'} \frac{\|K\|^2 |I|^2 x^3}{A^2_K n^2 \nu^2 T^2 h_{\min}} + \frac{ \|K\|^2 |I|^2 x^3}{2 h_{\min} A^2_K n^2 \nu^2 T^2} + \frac{ \|K\|_1^2 \|q\|_{I,\infty} \|l_{T}\|_{I,\infty} |I|^2}{n A^2_K \nu^2 T^2} \\
&\leq \frac{\theta'}{3} V_h  + \frac{1}{ \theta'} \frac{\|K\|_{\infty} \|K\|_1 |I|^2 x^3}{A^2_K n^2 \nu^2 T^2 h_{\min}} +  \frac{ \|K\|_1^2 \|q\|_{I,\infty} \|l_T\|_{I,\infty}|I|^2}{ \theta' n A^2_K \nu^2 T^2}. 
\end{align*}
To dominate 
\[C^2 = \frac{1}{n^4}\int_{0}^T \int_{0}^s \left(G_{h,h_{\min}}\left(X_u,X_s\right)+G_{h,h_{\min}}\left(X_s,X_u\right)\right)^2 d\Lambda_u d\Lambda_s\] which is bounded by 
\begin{equation} \label{boundCU}
\frac{4}{n^4} \int_{0}^T \int_{0}^T  \left(G_{h,h_{\min}}\left(X_u,X_s\right) \right)^2 d\Lambda_s d\Lambda_u,
\end{equation}
we use the occupation time formula and Young's inequality for convolutions:
\begin{align*}
C^2 &\leq \frac{4 \|q\|_{I,\infty} \|l_{T}\|_{I,\infty}|I|^2}{n^3 A^2_K \nu^2 T^2} \int_{0}^T \| w\left(\cdot,h,\frac{X_u-\cdot}{h}\right) K_h\left(X_u - \cdot \right){\bf 1}_{X_u \in I}{\bf 1}_{\cdot \in I} \ast K_{h_{\min}}\left(\cdot\right)  \|^2 d\Lambda_u\\
&\leq  \frac{4 \|q\|_{I,\infty} \|l_{T}\|_{I,\infty} \|K\|^2_1 |I|^2}{n^3 A^2_K \nu^2 T^2} \int_{0}^T \|w\left(\cdot,h,\frac{X_u-\cdot}{h}\right) K_h\left(X_u-\cdot \right){\bf 1}_{\cdot \in I}{\bf 1}_{X_u \in I}\|^2 d\Lambda_u\\
&=  \frac{4 \|q\|_{I,\infty} \|l_{T}\|_{I,\infty} \|K\|^2_1 |I|^2}{A^2_K n \nu^2 T^2} V_h.
\end{align*}
We then have
\[C\sqrt{x} \leq \frac{\theta'}{3} V_h + \tilde{C} \frac{\|K\|_1^2 \|q\|_{I,\infty} \|l_{T}\|_{I,\infty} |I|^2 x}{n A^2_K \nu^2 T^2 \theta'}.\]
Using twice Cauchy Schwarz inequality, we find  that 
\[D = \underset{\int_{0}^T a^2_u d\Lambda_u = 1, \int_{0}^T b^2_s d\Lambda_s = 1}{\sup} \int_{0}^T a_u \int_{u}^T b_s  \left(G_{h,h_{\min}}\left(X_u,X_s\right)+G_{h,h_{\min}}\left(X_s,X_u\right)\right) d\Lambda_s d\Lambda_u \]
is bounded by \eqref{boundCU} and then 
\[Dx \leq \frac{\theta'}{3} V_h + \tilde{C} \frac{\|K\|_1^2 \|q\|_{I,\infty} \|l_{T}\|_{I,\infty} |I|^2 x^2}{n A^2_K \nu^2 T^2 \theta'}.\]
Finally, with conditional probability given $\mathcal{F}^X_T$ larger than $1-\tilde{C} |\mathcal{H}|e^{-x}$, on $D\left(I,\nu\right)$,
\begin{equation} \label{controlU}
|U\left(h,h_{\min}\right)| \leq \theta' V_h + \tilde{C} \frac{\|K\|^2_1 \|q\|_{I,\infty} \|l_{T}\|_{I,\infty} |I|^2 x^2}{n  \nu^2 T^2 \theta'} + \tilde{C} \frac{\|K\|_{\infty}\|K\|_1 |I|^2 x^3}{\nu^2 T^2 A^2_K \theta' n^2 h_{\min}}.
\end{equation}

\paragraph{{\bf Control of $S$}} We need to control $S\left(h,h_{\min}\right)$ and $S\left(h_{\min},h\right)$. Let $h, h'$ in $\mathcal{H}$. We can write 
\[S\left(h,h'\right) = \frac{1}{n} \int_{0}^T \int_{I} \left(q_{h'}-q\right) w\left(x,h,\frac{X_s-x}{h}\right)K_{h}\left(X_s-x\right){\bf 1}_{X_s \in I} dx dM_s\]
and it is possible to control it with Proposition \ref{bernstein}. First, using occupation time formula and noticing that
\[q\left(x\right) - q_h\left(x\right) = \int_0^T \left(q\left(x\right)-q\left(X_s\right)\right)w\left(x,h,\frac{X_s-x}{h}\right) K_{h}\left(X_s-x\right){\bf 1}_{X_s \in I} ds  
\]
using Proposition \ref{polynom}, we have
\begin{align*}
\frac{1}{n}\int_{I} \left(q_{h'}-q\right)\left(x\right)  w\left(x,h,\frac{X_s-x}{h}\right) K_{h}\left(X_s-x\right){\bf 1}_{X_s \in I}dx &\leq \frac{|I|}{n A_K \nu T}\|K\|_1 \|q-q_{h'}\|_{I,\infty}\\
&\leq \frac{2 \|K\|^2_1 \|l_{T}\|_{I,\infty} \|q\|_{I,\infty} |I|^2}{n A^2_K \nu^2 T^2}. \\
\end{align*}
Using the occupation time formula and Young's inequality for convolutions, the term 
\[
\frac{1}{n^2}\int_{0}^T \left(\int_{I} \left(q_{h'}-q\right)\left(x\right)  w\left(x,h,\frac{X_s-x}{h}\right) K_{h}\left(X_s-x\right){\bf 1}_{X_s \in I} dx\right)^2d\Lambda_s\]
is bounded by 
\[ \frac{ \|l_{T}\|_{I,\infty} \|q\|_{I,\infty} |I|^2}{n A^2_K \nu^2 T^2} \| \left(q_{h'}-q\right)\left(\cdot\right){\bf 1}_{\cdot \in I} \ast K_h\left(\cdot\right) \|^2 \leq \frac{ \|q_{h'}-q\|_I^2 \|K\|_1^2 \|l_{T}\|_{I,\infty} \|q\|_{I,\infty} |I|^2}{n A^2_K \nu^2 T^2}. \]
With conditional probability given $\mathcal{F}^X_T$ larger than $1-\tilde{C} e^{-x}$, on $D\left(I,\nu\right)$, we then have
\begin{equation} \label{controlS}
S\left(h,h'\right) \leq \frac{\theta'}{2}\|q_{h'}-q\|_I^2 + \tilde{C} \frac{\|K\|_1^2 \|l_{T}\|_{I,\infty} |I|^2 \|q\|_{I,\infty} x}{n A^2_K \nu^2 T^2 \theta'}.
\end{equation}
We apply \eqref{controlS} for $S\left(h,h_{\min}\right)$ and $S\left(h_{\min},h\right)$.

\medskip
\paragraph{{ \bf Control of $<q_{h}-q,q_{h_{\min}}-q>_I$}} We have: 
\begin{equation} \label{controlscalar}
| <q_{h}-q,q_{h_{\min}}-q>_I| \leq \frac{\theta'}{2}\|q_{h}-q\|_I^2 + \frac{1}{2\theta'}\|q_{h_{\min}}-q\|_I^2.
\end{equation}

At the end, combining \eqref{decompositionpen}, \eqref{controlU}, \eqref{controlS} and \eqref{controlscalar}, we have with conditional probability given $\mathcal{F}^X_T$ higher than $1-\tilde{C} |\mathcal{H}| e^{-x}$, on $D\left(I,\nu\right)$, for any $h \in \mathcal{H}$,
\begin{equation} 
\label{controlpen2}
\begin{split}
|<\hat{q}_{h}-q,\hat{q}_{h_{\min}}-q>_I - \hat{V}_{h,h_{\min}}| &\leq \theta'\left(\|q_{h}-q\|_I^2+V_h\right)+ \left(\frac{\theta'}{2}+\frac{1}{2\theta'}\right)\|q_{h_{\min}}-q\|_I^2 \\
&+ \frac{\tilde{C}\left(K\right)|I|}{\nu^2 T^2  \theta'}\left(\frac{\|q\|_{I,\infty} \|l_{T}\|_{I,\infty} |I| x^2}{n} + \frac{x^3 |I|}{n^2 h_{\min}}\right). 
\end{split}
\end{equation}
Furthermore, Proposition \ref{biasvariance} states that with conditional probability given $\mathcal{F}^X_T$ higher than $1-\tilde{C} |\mathcal{H}| e^{-x}$ for any $h \in \mathcal{H}$, on $D\left(I,\nu\right)$,
\begin{equation}
\|q- q_{h}\|_I^2 + V_h  \leq 2 \|q- \hat{q}_{h}\|_I^2 + \tilde{C} \frac{|I|\left(1+ \|q\|_{I,\infty}\|l_T\|_{I,\infty}|I| \|K\|_1^2\right)x^2}{nA_K^2\nu^2T^2}. \label{controlpen3}
\end{equation}
We combine \eqref{controlpen}, \eqref{controlpen2} and \eqref{controlpen3} and we take $\theta' = \frac{\theta}{4}$ to conclude.

\subsection{Proof of Proposition \ref{oracle}}

\paragraph{{\bf Proof of \eqref{oracleinequality}}} Let $x \geq 1$, $\tau = \alpha-1$, $\epsilon \in \left(0,1\right)$ and $\theta \in \left(0,1\right)$ depending on $\epsilon$ and specified later. Using Proposition \ref{oraclewithpen}, with conditional probability given $\mathcal{F}^X_T$ larger than $1-\tilde{C} |\mathcal{H}| e^{-x}$, on $D\left(I,\nu\right)$, for any $h \in \mathcal{H}$,
\begin{equation} \label{oracleinequalityprep1} \begin{split}
\left(1-\theta\right)\|\hat{q}_{\hat{h}}-q\|_I^2 + \tau \hat{V}_{\hat{h}} \leq &\left(1+\theta\right)\|\hat{q}_{h}-q\|_I^2 + \tau \hat{V}_h + \frac{\tilde{C}}{\theta}\|q_{h_{\min}}-q\|_I^2 \\
&+\frac{\tilde{C}\left(K \right)|I|}{\nu^2 T^2 \theta}\left(\frac{\left(1+\|q\|_{I,\infty} \|l_T\|_{I,\infty}|I|\right) x^2}{n} + \frac{x^3|I|}{n^2 h_{\min}}   \right)
\end{split}
\end{equation}
Equation \eqref{variancepart1} states that with conditional probability given $\mathcal{F}^X_T$ larger than $1-\tilde{C}|\mathcal{H}| e^{-x}$, on $D\left(I,\nu\right)$, for any $h \in \mathcal{H}$,
\begin{equation} \label{oracleinequalityprep2}
|\hat{V}_h - V_h | \leq \frac{\theta}{2} V_h +  \tilde{C} \frac{x|I|}{n A_K^2 \nu^2 T^2 \theta}.
\end{equation}
First, let us consider the case $\tau \geq 0$. We then have, with conditional probability given $\mathcal{F}^X_T$ larger than $1-\tilde{C}|\mathcal{H}| e^{-x}$, on $D\left(I,\nu\right)$, for any $h \in \mathcal{H}$,
\[
\begin{split}
 \left(1-\theta\right)\|\hat{q}_{\hat{h}}-q\|_I^2  &\leq \left(1+\theta\right)\|\hat{q}_{h}-q\|_I^2 + \tau \left(1+\frac{\theta}{2}\right)V_h + \frac{\tilde{C}}{\theta}\|q_{h_{\min}}-q\|_I^2 \\
&+\frac{\tilde{C}\left(K \right)|I|}{\nu^2 T^2 \theta}\left(\frac{\left(1+\|q\|_{I,\infty} \|l_T\|_{I,\infty}|I|\right) x^2 }{n} + \frac{|I|x^3}{n^2 h_{\min}}   \right).
\end{split}
\]
Using Proposition \ref{biasvariance} with $\eta = \frac{\theta}{2+\theta}$, with conditional probability given $\mathcal{F}^X_T$ larger than $1-\tilde{C}|\mathcal{H}| e^{-x}$, on $D\left(I,\nu\right)$, for any $h \in \mathcal{H}$,
\begin{equation} \label{oracleinequalityprep3}
\tau V_h \leq \tau\left(1+\frac{\theta}{2+\theta}\right)\|q-\hat{q}_h\|_I^2 + \tau \frac{\tilde{C}\left(K\right)|I|\left(1+\|q\|_{I,\infty}\|l_T\|_{I,\infty}|I|\right) x^2}{n \theta^3 \nu^2 T^2}.
\end{equation}
Combining inequality \eqref{oracleinequalityprep3} with \eqref{oracleinequalityprep1} and \eqref{oracleinequalityprep2} and as $\left(1+\frac{\theta}{2}\right)\left(1+\frac{\theta}{2+\theta}\right) = 1+\theta$, we find that with conditional probability given $\mathcal{F}^X_T$ larger than $1-\tilde{C}|\mathcal{H}| e^{-x}$, on $D\left(I,\nu\right)$, for any $h \in \mathcal{H}$,
\[
\begin{split}
  \left(1-\theta\right)\|\hat{q}_{\hat{h}}-q\|_I^2  &\leq \left(1+\theta + \left(1+\theta\right)\tau \right)\|\hat{q}_{h}-q\|_I^2 + \frac{\tilde{C}}{\theta}\|q_{h_{\min}}-q\|_I^2 \\
&+\frac{|I|\tilde{C}\left(K \right)}{\nu^2 T^2}\left(\frac{1}{\theta}+\frac{1}{\theta^3}\right)\left(\frac{\left(1+\|q\|_{I,\infty}\|l_T\|_{I,\infty}|I|\right)x^2}{n} + \frac{|I|x^3}{n^2 h_{\min}}\right).
\end{split}
\]
With $\theta = \frac{\epsilon}{\epsilon+2+2\tau}$, we have
\[
\begin{split} \|\hat{q}_{\hat{h}}-q\|_I^2 \leq &\left(1+\tau+\epsilon\right) \|\hat{q}_h-q\|_I^2 +  \frac{\tilde{C}\left(\epsilon+2+2\tau\right)^2}{\left(2+2\tau\right)\epsilon}\|q_{h_{\min}}-q\|_I^2\\
& + \frac{\tilde{C}\left(K\right)|I|\left(\epsilon+2+2\tau\right)^4}{\nu^2 T^2\left(2+2\tau\right)\epsilon^3}\left(\frac{\left(1+\|q\|_{I,\infty}\|l_T\|_{I,\infty}|I|\right)x^2}{n} + \frac{x^3|I|}{n^2 h_{\min}}\right).
\end{split}
\]
Inequality \eqref{oracleinequality} is then proved in the case $\alpha \geq 1$ with
\[C_0\left(\epsilon,\alpha\right) = \alpha + \epsilon,\]
\[C_2\left(\epsilon,\alpha\right) = \tilde{C} \frac{\left(\epsilon+2\alpha\right)^2}{\alpha \epsilon}\]
and
\[C_3\left(\epsilon,K,\alpha\right) = \tilde{C}\left(K\right) \frac{\left(\epsilon+2\alpha\right)^4}{\alpha \epsilon^3}.\]

\medskip
Now, let us consider the case $-1 < \tau \leq 0$. According to Proposition \ref{biasvariance}, with conditional probability on $\mathcal{F}^X_T$ larger than $1-\tilde{C}|\mathcal{H}| e^{-x}$, on $D\left(I,\nu\right)$, 
\begin{equation} \label{oracleinequalityprep4}
\tau V_{\hat{h}} \geq \tau\left(1+\theta\right)\|q-\hat{q}_{\hat{h}}\|_I^2 + \tau \frac{\tilde{C}\left(K\right)|I|\left(1+\|q\|_{I,\infty}\|l_T\|_{I,\infty}|I|\right) x^2}{n \theta^3 \nu^2 T^2}.
\end{equation}
We find, combining \eqref{oracleinequalityprep4} with \eqref{oracleinequalityprep1} and \eqref{oracleinequalityprep2} , that with conditional probability given $\mathcal{F}^X_T$ larger than $1-\tilde{C}|\mathcal{H}| e^{-x}$, on $D\left(I,\nu\right)$, for any $h \in \mathcal{H}$,
\[
\begin{split}
\left(1-\theta + \tau\left(1+\theta\right)\right) &\|\hat{q}_{\hat{h}}-q\|_I^2  \leq \left(1+\theta \right)\|\hat{q}_{h}-q\|_I^2 + \frac{\tilde{C}}{\theta}\|q_{h_{\min}}-q\|_I^2 \\
&+\frac{|I|\tilde{C}\left(K\right)}{\nu^2 T^2}\left(\frac{1}{\theta}+\frac{1}{\theta^3}\right)\left(\frac{\left(1+\|q\|_{I,\infty}\|l_T\|_{I,\infty}|I|\right)x^2}{n} + \frac{x^3|I|}{n^2 h_{\min}}\right).
\end{split}
\]
Taking $\theta =\frac{\epsilon\left(\tau+1\right)^2}{2+\epsilon\left(1-\tau^2\right)} < 1$, we obtain, with conditional probability given $\mathcal{F}^X_T$ larger than $1-\tilde{C}|\mathcal{H}| e^{-x}$, on $D\left(I,\nu\right)$, for any $h \in \mathcal{H}$,
\[
\begin{split}
 \|\hat{q}_{\hat{h}}-q\|_I^2 \leq &\left(\frac{1}{1+\tau}+\epsilon\right) \|\hat{q}_h-q\|_I^2 +  \frac{\tilde{C}\left(2+\epsilon\left(1-\tau\right)\right)^2}{2\epsilon\left(\tau+1\right)^3}\|q_{h_{\min}}-q\|_I^2\\
& + \frac{\tilde{C}\left(K\right) \left(2+\epsilon\left(1-\tau\right)\right)^4 |I|}{2\epsilon^3\left(\tau+1\right)^7}\left(\frac{\left(1+\|q\|_{I,\infty}\|l_T\|_{I,\infty}|I|\right) x^2}{n} + \frac{x^3|I|}{n^2 h_{\min}}\right).
\end{split}
\]
Inequality \eqref{oracleinequality} is then verified in the case $\alpha < 1$ with
\[C_0\left(\epsilon,\alpha\right) = \frac{1}{\alpha} + \epsilon,\]
\[C_2\left(\epsilon,\alpha\right) = \tilde{C} \frac{\left(2+\epsilon\left(2-\alpha\right)\right)^2}{2\epsilon \alpha^3}\]
and
\[C_3\left(\epsilon,K,\alpha\right) = \tilde{C}\left(K\right)  \frac{\left(2+\epsilon\left(2-\alpha\right)\right)^4}{2\epsilon^3\alpha^7}.\]

\medskip
\paragraph{{\bf Proof of \eqref{oracleexpectation}}} Let us use \eqref{oracleinequality} with $x = 5\log\left(n \vee |\mathcal{H}|\right)$ and $\epsilon = \tilde{C} \left(\log\left(n\right)\right)^{-1}$. Let $\mathcal{E}$ be the event on which \eqref{oracleinequality} is true. Integrating with respect to $\mathcal{F}^X_T$ and dividing by $\mathbb{P}\left(D\left(I,\nu\right)\right)$, we find 
\begin{equation}\label{controlexpectation1}
\begin{split}
&\mathbb{E}\left(\|\hat{q}_{\hat{h}}-q\|_I^2 {\bf 1}_{\mathcal{E}} | D\left(I,\nu\right)\right) \leq \left(\alpha \vee \frac{1}{\alpha} +\tilde{C} \left(\log\left(n\right)\right)^{-1} \right) \underset{h \in \mathcal{H}}{\min} \; \mathbb{E}\left( \| \hat{q}_{h} - q\|_I^2| D\left(I,\nu\right)\right) \\
&+\tilde{C}\left(\alpha\right)\log\left(n\right)\mathbb{E}\left(\|q_{h_{\min}} - q\|_I^2 | D\left(I,\nu\right)\right)\\
&+ \frac{\tilde{C}\left(K,\alpha\right)|I|\log\left(n\right)^{3}}{\nu^2 T^2}\left(\frac{25|I| \|q\|_{I,\infty} \mathbb{E}\left(\|l_T\|_{I,\infty} | D\left(I,\nu\right)\right) \log\left(n \vee |\mathcal{H}|\right)^2}{n}  +\frac{125|I| \log\left(n \vee |\mathcal{H}|\right)^3}{n^2 h_{\min}}\right).    
\end{split}
\end{equation}
Let $N_I = \int_{0}^T {\bf 1}_{X_s \in I} dN_s$. On $\mathcal{E}^c \cap D\left(I,\nu\right)$, using Cauchy Schwarz inequality,
\begin{align*}
\|\hat{q}_{\hat{h}}-q\|_I^2 &\leq 2 \| q \|^2_I + \frac{2|I|^2}{\nu^2 T^2 A^2_K n^2} N_I \int_{I}  \int_{0}^T K^2_h\left(X_s-x\right){\bf 1 }_{X_s \in I} dN_s dx\\
&\leq 2 \| q \|^2_I + \frac{2 N_I^2 \|K\|_{\infty} \|K\|_1|I|^2}{\nu^2 T^2 A^2_K n^2 h}\\
&\leq 2 \| q \|^2_I + \frac{2 N_I^2 |I|}{\nu^2 T^2 A^2_K n}.
\end{align*}
Using Cauchy Schwarz inequality again, we have
\[\mathbb{E}\left(N_I^2{\bf 1}_{\mathcal{E}^c} | \mathcal{F}^X_T \right) \leq \mathbb{E}\left(N_I^4 | \mathcal{F}^X_T \right)^{\frac{1}{2}} \mathbb{P}\left(\mathcal{E}^c | \mathcal{F}^X_T\right)^{\frac{1}{2}}.\] 
Using Laplace transform formula, we easily show that $N_I$ has the law of a Poisson random variable with parameter $n \int_0^T q\left(X_s\right) {\bf 1}_{X_s \in I} ds$ conditionally on $\mathcal{F}^X_T$ and 
\[\mathbb{E}\left(N_I^4 | \mathcal{F}^X_T\right) \leq \tilde{C} n^4 \left( \|q\|_{I,\infty} T+\left(\|q\|_{I,\infty} T\right)^2 +\left(\|q\|_{I,\infty} T\right)^3 +\left(\|q\|_{I,\infty} T\right)^4\right).\]
We also have 
\begin{align*} \mathbb{P}\left(\mathcal{E}^c | \mathcal{F}^X_T\right) &\leq \tilde{C} \frac{|\mathcal{H}|}{n^5 \vee |\mathcal{H}|^5}\\
&\leq \frac{\tilde{C}}{n^4}.
\end{align*}
Integrating with respect to $\mathcal{F}^X_T$ and dividing by $\mathbb{P}\left(D\left(I,\nu\right)\right)$, we find 
\begin{equation} \label{controlexpectation2}
\mathbb{E}\left(\|\hat{q}_{\hat{h}}-q\|_I^2 {\bf 1}_{\mathcal{E}^c} | D\left(I,\nu\right)\right) \leq \tilde{C}\frac{\|q\|^2_I}{n^4} + \tilde{C}\left(K\right)\frac{|I|\sqrt{\sum_{i=1}^4 \left(\|q\|_{\infty,I}T\right)^i}}{\nu^2 T^2 n}. 
\end{equation}
We obtain \eqref{oracleexpectation} combining \eqref{controlexpectation1} and \eqref{controlexpectation2}.

\section{Other proofs}
\label{otherproofs}

\subsection{Proof of Proposition \ref{propertiesX}}

Under Assumption \ref{assumptionsX} (i), we can define a local time in the sense of the continuous semimartingale $L_t^x$ continuous in $t$, cadlag in $x$, see \cite[Chapter 6]{revuz13} for more information. \cite[Exercise 1.15]{revuz13} states that for every measurable function $h$ on $\left[0,T\right] \times \Omega \times \mathbb{R}$. 
\[\int_0^t h\left(s,X_s\right)d<X_s,X_s> = \int_{\mathbb{R}} da \int_{0}^t h\left(s,a\right)dL_s^a.\]
Let $f$ be a measurable function on $\Omega \times \mathbb{R}$. As $\sigma_s \geq \underline{\sigma} > 0$  for every $s \in \left[0,T\right]$ almost surely and $d<X_s,X_s> = \sigma_s^2 ds$, we have 
\[\int_0^t f\left(X_s\right)ds = \int_{\mathbb{R}} f\left(a\right) da \int_{0}^t \frac{1}{\sigma_s^2} dL_s^a\]
and (i) is verified with $l_T^x = \int_{0}^T \frac{1}{\sigma_s^2} dL_s^a$. According to \cite[Equation $(\text{III})_{\gamma}$]{barlow82}, as 
\[\mathbb{E}\left(\underset{0 \leq t \leq T}{\sup}|\int_{0}^t \sigma_s dW_s| + \int_{0}^T |\mu_s| ds\right) < \infty,\] 
\[\mathbb{E}\left(\underset{x \in \mathbb{R}}{\sup}\; L_T^x\right) < \infty\]
and because $\sigma_s \geq \underline{\sigma}$ a.s., we obtain (ii). The continuity of $x \mapsto l_T^x$ follows from \cite[Example 2.2.3 (a)]{yen13} and $\sigma^2_s > 0$ for all $s \in \left[0,T\right]$.

\medskip
Under Assumption \ref{assumptionsX} (ii), we have $l_t^x = {\bf 1}_{x \in \left[0,t\right]}$ for all $x$ in $\mathbb{R}$ and $t$ in $\left[0,T\right]$.

\subsection{Proof of Proposition \ref{minB}}

Let $\|\cdot\|_{m+1}$ be the Euclidian norm on $\mathbb{R}^{m+1}$. For a symmetric positive semi-definite matrix $A$, let $\lambda_{\min}\left(A\right)$ the smallest eigenvalue of $A$. In the following, we work on $D\left(I,\nu\right)$. We have
\[\lambda_{\min}\left(B\left(x,h\right)\right) = \underset{ \|v\|_{m+1} = 1}{\inf} v^T B\left(x,h\right) v.\]
Let $v$ in $\mathbb{R}^{m+1}$ with $\|v\|_{m+1} = 1$. Using the occupation time formula, we have
\begin{align*}
v^T B\left(x,h\right) v &=  \int_0^T \left(v^T U\left(\frac{X_s-x}{h}\right)\right)^2 K_h\left(X_t-x\right){\bf 1}_{X_s \in I}  ds \\
&= \int_{I} \left(v^T U\left(\frac{u-x}{h}\right)\right)^2 K_h\left(u-x\right)  l_T^u du\\
&\geq \frac{\nu T K_{\min}}{|I|h} \int_{I} \left(v^T U\left(\frac{u-x}{h}\right)\right)^2 {\bf 1}_{|u-x| < \Delta h} du.
\end{align*}
If $x \in \left[\min I + \Delta h , \max I - \Delta h\right]$,
\[\frac{1}{h}\int_{I} \left(v^T U\left(\frac{u-x}{h}\right)\right)^2 {\bf 1}_{|u-x| < \Delta h} du = \int_{\mathbb{R}} \left(v^T U\left(u\right)\right)^2 {\bf 1}_{|u| < \Delta}du ,\]
if $x \in \left[\min I , \min I  + \Delta h\right]$, because $h \leq \frac{2}{3}  \frac{| I |}{\Delta}$,
\[\frac{1}{h}\int_{I} \left(v^T U\left(\frac{u-x}{h}\right)\right)^2 {\bf 1}_{|u-x| < \Delta h}du \geq \int_{\mathbb{R}} \left(v^T U\left(u\right)\right)^2 {\bf 1}_{0 < u < \frac{\Delta}{2}}du \]
and if $x \in \left[\max I -  \Delta h , \max I\right]$,
\[\frac{1}{h}\int_{I} \left(v^T U\left(\frac{u-x}{h}\right)\right)^2 {\bf 1}_{|u-x| < \Delta h} du\geq \int_{\mathbb{R}} \left(v^T U\left(u\right)\right)^2 {\bf 1}_{\frac{-\Delta}{2} < u < 0}du.\]
Thus, for all $x \in I$, $\lambda_{\min}\left(B\left(x,h\right)\right)$ is larger than 
\begin{align*}
 &\frac{\nu T}{|I|} K_{\min} \min \left(\underset{ \|v\|_{m+1} = 1}{\inf} \left(\int_{\mathbb{R}} \left(v^T U\left(u\right) \right)^2 {\bf 1}_{0 < u < \frac{\Delta}{2}} du\right),\underset{ \|v\|_{m+1} = 1}{\inf} \left(\int_{\mathbb{R}} \left(v^T U\left(u\right) \right)^2 {\bf 1}_{\frac{-\Delta}{2} < u < 0} du\right)\right)  \\
&= \frac{\nu T}{|I|} K_{\min} \min\left( \lambda_{\min}\left( \int_{\mathbb{R}} U\left(u\right) U^T\left(u\right) {\bf 1}_{0 < u < \frac{\Delta}{2}} du \right),\lambda_{\min}\left( \int_{\mathbb{R}} U\left(u\right) U^T\left(u\right) {\bf 1}_{\frac{-\Delta}{2} < u < 0} du \right)\right) 
\end{align*}
and 
\[ \min\left( \lambda_{\min}\left( \int_{\mathbb{R}} U\left(u\right) U^T\left(u\right) {\bf 1}_{0 < u < \frac{\Delta}{2}} du \right),\lambda_{\min}\left( \int_{\mathbb{R}} U\left(u\right) U^T\left(u\right) {\bf 1}_{\frac{-\Delta}{2} < u < 0} du \right)\right) > 0\]
applying \cite[Lemma 1.4]{tsybakov09} with $K\left(u\right) = {\bf 1}_{0 \leq u \leq \frac{\Delta}{2}}$ and $K\left(u\right) = {\bf 1}_{-\frac{\Delta}{2} \leq u \leq 0}$.

We have 
\begin{align*}
|w\left(x,h,\frac{X_s-x}{h}\right) {\bf 1}_{|\frac{X_s-x}{h}| \leq 1}| &\leq \|B\left(x,h\right)^{-1} U\left(\frac{X_s-x}{h}\right)  \|_{m+1} {\bf 1}_{|\frac{X_s-x}{h}| \leq 1}\\
&\leq  \frac{1}{\lambda_{\min}\left(B\left(x,h\right)\right)} \|U\left(\frac{X_s-x}{h}\right) \|_{m+1}{\bf 1}_{|\frac{X_s-x}{h}| \leq 1}\\
&\leq \frac{1}{\lambda_{\min}\left(B\left(x,h\right)\right)} \sqrt{1 + \frac{1}{\left(1!\right)^2} + \frac{1}{\left(2!\right)^2}+...+\frac{1}{\left(m!\right)^2}}\\
&\leq \frac{2}{\lambda_{\min}\left(B\left(x,h\right)\right)}.
\end{align*}
Then,
\begin{equation} \label{boundw}
 |w\left(x,h,\frac{X_s-x}{h}\right) {\bf 1}_{|\frac{X_s-x}{h}| \leq 1}| \leq \frac{|I|}{A_K \nu T}
\end{equation}
with 
\[A_K = \frac{ K_{\min} \min\left( \lambda_{\min}\left( \int_{\mathbb{R}} U\left(u\right) U^T\left(u\right) {\bf 1}_{0 < u < \frac{\Delta}{2}} du \right),\lambda_{\min}\left( \int_{\mathbb{R}} U\left(u\right) U^T\left(u\right) {\bf 1}_{-\frac{\Delta}{2} < u < 0} du \right)\right)}{2} > 0.\]

\subsection{Proof of Proposition \ref{minimaxnonparametric}}

\paragraph{{\bf Lower bound.}} We suppose for simplicity that $I = \left[0,1\right]$. As in \cite[Section 2.6.1]{tsybakov09}, we consider a real number $c_0 > 0$ and 
\begin{align*}
&m = \lfloor c_0 n^{\frac{1}{2\beta+1}}\rfloor + 1,\; h_n = \frac{1}{m}, \; x_k = \frac{k-\frac{1}{2}}{m},\\
&\varphi_k\left(x\right) = Lh_n^{\beta}\overline{K}\left(\frac{x-x_k}{h_n}\right), \; k = 1,..,m,\; x \in \left[0,1\right]
\end{align*}
with 
\[\overline{K}\left(u\right) = a \exp\left(-\frac{1}{1-4u^2}\right){\bf 1}_{|2u| \leq 1}, \; a > 0.\] 
For $a$ sufficient small, $\overline{K} \in \Sigma\left(\beta, \frac{1}{2} \right) \cap \mathcal{C}^{\infty}\left(\mathbb{R}\right)$. According to \cite[Equation (2.5)]{tsybakov09}, the functions $\varphi_k$ belongs to $\Sigma\left(\beta,\frac{L}{2},I\right)$ and the set of function 
\[\mathcal{C} = \{q: \; q\left(x\right) = \rho +  \sum_{k=1}^{m} w_k \varphi_{k}\left(x\right), \; w_k \in  \{0,1\},\; x \in I\}\]
is included in $\Lambda_{\rho,\beta}$  as the functions $\varphi_k$ have disjoint supports.

\medskip
Let us suppose that $m \geq 8$. According to \cite[Lemma 2.9]{tsybakov09}, there exists a subset $\mathcal{\tilde{C}}$ of $\mathcal{C}$ such that for all $f_w  =  \rho + \sum_{k=1}^{m} w_k \varphi_{k}\left(x\right) \in \mathcal{\tilde{C}}$ and all  $f_{w'} = \rho+ \sum_{k=1}^{m} w'_k \varphi_{k}\left(x\right) \in \mathcal{\tilde{C}}$, we have 
\[\sum_{k=1}^{m} \left(w_k-w'_k\right)^2 \geq \frac{m}{8}\]
and with 
\[M \geq \frac{2^{m}}{8}\]
where $M+1 = |\mathcal{\tilde{C}} |$.
Now, if we consider two elements $f_w$ and $f_{w'}$ of $\mathcal{\tilde{C}}$, we have
\begin{align*}
\|f_w-f_{w'}\|_I &=L h_n^{\beta + \frac{1}{2}} \|\overline{K}\| \sqrt{\sum_{k=1}^{m} \left(w_k-w'_k\right)^2}\\
&\geq L h_n^{\beta + \frac{1}{2}} \|\overline{K}\| \sqrt{\frac{m}{16}}\\
& = \frac{L}{4} \|\overline{K}\| m^{-\beta}.
\end{align*}
Thus, if $n \geq n^* =  \left(\frac{7}{c_0}\right)^{2\beta+1}$, $m \geq 8$ and $m^{\beta} \leq \left(2c_0\right)^{\beta } n^{\frac{\beta}{2\beta+1}}$. Hence, 
\begin{equation} \label{conddistance}
\|f_w-f_{w'}\|_I \geq 2s_n
\end{equation}
with $s_n = A n^{-\frac{\beta}{2\beta+1}}$ and $A = \frac{L}{8}\|\overline{K}\|\left(2c_0\right)^{-\beta}$.

\medskip
The following part of the proof differs from \cite{tsybakov09}. $\mathcal{\tilde{C}}$ can be written $\{q_0, q_1, ...q_M\}$. Let us denote by $\mathbb{P}_j$ the probability measure associated to the intensity $n q_j\left(X_s\right), \; s \in \left[0,T\right]$. We consider the Kullback divergence between $\mathbb{P}_0$ and $\mathbb{P}_j$ conditionally on $D\left(I,\nu\right)$, that is $\mathbb{E}^{\mathbb{P}_0}\left(\log\left(\frac{d\mathbb{P}_j}{d\mathbb{P}_0}\right) | D\left(I,\nu\right)\right)$, denoted by $K\left(\mathbb{P}_{0},\mathbb{P}_j\right)$. We have
\begin{align*}
K\left(\mathbb{P}_0, \mathbb{P}_j\right) &= \mathbb{E}\left(\int_{0}^T n \left(q_j\left(X_s\right)-q_0\left(X_s\right)- q_0\left(X_s\right) \log\left(\frac{q_j\left(X_s\right)}{q_0\left(X_s\right)} \right)\right)ds | D\left(I,\nu\right) \right) \\
&=\mathbb{E}\left(\int_{0}^T n q_0\left(X_s\right)\left(\frac{q_j\left(X_s\right)-q_0\left(X_s\right)}{q_0\left(X_s\right)} - \log\left(1+\frac{q_j\left(X_s\right)-q_0\left(X_s\right)}{q_0\left(X_s\right)} \right)\right)ds | D\left(I,\nu\right) \right)  \\
&\leq \mathbb{E}\left(\int_{0}^T n q_0\left(X_s\right)\left(\frac{q_j\left(X_s\right)-q_0\left(X_s\right)}{q_0\left(X_s\right)} - \log\left(1+\frac{q_j\left(X_s\right)-q_0\left(X_s\right)}{q_0\left(X_s\right)} \right)\right)ds | D\left(I,\nu\right) \right)  \\
&\leq \mathbb{E}\left(\int_{0}^T n \frac{\left(q_j\left(X_s\right)-q_0\left(X_s\right)\right)^2}{q_j\left(X_s\right)} ds | D\left(I,\nu\right) \right) \numberthis \label{kullbackbounding}
\end{align*}
using the fact that for $x > -1$, $\log\left(1+x\right) \geq \frac{x}{1+x}$. 
Continuing from \eqref{kullbackbounding}, we have:
\begin{align*}
K\left(\mathbb{P}_0, \mathbb{P}_j\right) &\leq  \frac{n h_n^{2\beta} L^2 T \|\overline{K}\|^2_{\infty}}{\rho}  \\
&=\frac{L^2 T \|\overline{K}\|^2_{\infty} c_0^{-\left(2\beta+1\right)}m}{\rho}.  \\
\end{align*}
Let $\alpha \in \left(0,\frac{1}{8}\right)$. As $\log\left(M\right) \geq \frac{\log\left(2\right)m}{8}$, we choose 
\[c_0 = \left(\frac{8L^2\|\overline{K}\|^2_{\infty}T}{\alpha \log\left(2\right)\rho}\right)^{\frac{1}{2\beta + 1}}\]
and we have 
\begin{equation} \label{condkullback}
\frac{1}{M} \sum_{j=1}^M K\left(\mathbb{P}_0,\mathbb{P}_j\right) \leq \alpha \log\left(M\right).
\end{equation}
We have according to \eqref{conddistance} and \eqref{condkullback}
\begin{enumerate}
\item[(i)] $\|\lambda_j-\lambda_k\|_I > 2s_n$ for all $j \neq k$ with $s_n = A n^{-\frac{\beta}{2\beta+1}}$ and 
\item[(ii)] $\frac{1}{M} \sum_{j=1}^M K\left(\mathbb{P}_0,\mathbb{P}_j\right) \leq \alpha \log\left(M\right)$.
\end{enumerate}
We can conclude using \cite[Theorem 2.5]{tsybakov09}:
\[
\underset{n \to \infty}{\lim \inf} \; \underset{T_n}{\inf} \; {\bf R}\left(T_n, \Lambda_{\rho,\beta}, n^{\frac{-\beta}{2\beta+1}}\right) > 0.
\]

\paragraph{{\bf Upper bound.}} Let us assume that $l \leq m$. Let $h$ in $\mathcal{H}$. If $q \in \Sigma\left(L,\beta,I\right)$, the bias part of $\mathbb{E}\left(\|\hat{q}_h-q\|_I^2 | D\left(I,\nu\right)\right)$,  $\mathbb{E}\left(\|q-q_h\|_I^2 | D\left(I,\nu\right)\right) $, is equal to, using Proposition \ref{polynom}, 
\begin{equation} \label{biasequal}
\mathbb{E}\left( \int_I \left(\int_0^T  w\left(x,h,\frac{X_s-x}{h}\right) K_h\left(X_s-x\right){\bf 1}_{X_s \in I} \left(q\left(X_s\right)-q\left(x\right)\right) ds  \right)^2dx | D\left(I,\nu\right) \right).
\end{equation}
Using Proposition \ref{polynom} and Taylor's expansion, there exists $\left(\tau_s\right)_{0 \leq s \leq T}$ such that the integral inside the expectation in \eqref{biasequal} is equal to, on the event $D\left(I,\nu\right)$, 
\[
 \int_I \left(\int_0^T  w\left(x,h,\frac{X_s-x}{h}\right) K_h\left(X_s-x\right){\bf 1}_{X_s \in I} \left(q^{\left(l\right)}\left(x+\tau_s\left(X_s-x\right)\right)-q^{\left(l\right)}\left(x\right)\right)\frac{\left(X_s-x\right)^{l}}{l!} ds\right)^2dx. 
\]  
The bias term is then bounded by 
\begin{align*}
&\mathbb{E}\left( \int_I \left(\int_0^T  |w\left(x,h,\frac{X_s-x}{h}\right)| K_h\left(X_s-x\right){\bf 1}_{X_s \in I} \frac{L |X_s-x|^{\beta}}{l!}ds \right)^2dx | D\left(I,\nu\right) \right) \\
&\leq \mathbb{E}\left( \int_I \left(\int_0^T  |w\left(x,h,\frac{X_s-x}{h}\right)| K_h\left(X_s-x\right){\bf 1}_{X_s \in I} \frac{L h^{\beta}}{l!} ds\right)^2dx | D\left(I,\nu\right) \right) \\
&\leq \frac{|I|^2 L^2 h^{2\beta}}{\nu^2 T^2 A_K^2 \left(l!\right)^2} \mathbb{E}\left(\|l_T\|_{I,\infty}\int_I \int_I \int_0^T K_h\left(X_s-x\right)K_h\left(z-x\right){\bf 1}_{X_s \in I} ds l^z_T dz dx| D\left(I,\nu\right) \right)\\
&\leq \frac{ L^2 h^{2\beta}\mathbb{E}\left(\|l_T\|_{I,\infty} | D\left(I,\nu\right)\right)\|K\|_1^2 |I|^2}{\nu^2 T A_K^2 \left(l!\right)^2}.
\end{align*}
The variance part is equal to 
\begin{align*}
\mathbb{E}\left(\|\hat{q}_h-q_h\|_I^2 | D\left(I,\nu\right)\right) &= \frac{1}{n}\mathbb{E}\left(\int_I \int_0^T w^2\left(x,h,\frac{X_s-x}{h}\right) K^2_h\left(X_s-x\right) q\left(X_s\right) ds dx | D\left(I,\nu\right)  \right) \\
&\leq \frac{\|q\|_{I,\infty}|I|^2}{A^2_K \nu^2 T^2 n}\mathbb{E}\left(\int_I \int_0^T K^2_h\left(X_s-x\right) ds dx | D\left(I,\nu\right) \right)\\
&\leq \frac{\|q\|_{I,\infty} \|K\|^2 |I|^2}{nh A^2_K \nu^2 T}.
\end{align*}
Hence, we have
\[\underset{h \in \mathcal{H}}{\min} \; \mathbb{E}\left(\|\hat{q}_h-q \|^2_I | D\left(I,\nu\right)\right) \leq C_1\left(q,\mathbb{E}\left(\|l_T\|_{I,\infty} | D\left(I,\nu\right)\right),K,T,\beta,|I|\right) n^{\frac{-2\beta}{2\beta+1}} \]
and according to \eqref{oracleexpectation}, as $|\mathcal{H}| \leq \tilde{C}\left(K,|I|\right) n$ where $\tilde{C}\left(K,|I|\right)$ is a constant depending only on $K$ and $|I|$ and $\mathbb{E}\left(\|q-q_{h_{\min}}\|^2 | D\left(I,\nu\right)\right)$ is of order $n^{-2\beta}$, 
\[\mathbb{E}\left(\|\hat{q}_{\hat{h}}-q\|^2_I  | D\left(I,\nu\right) \right) \leq C_2\left(q,\mathbb{E}\left(\|l_T\|_{I,\infty} | D\left(I,\nu\right)\right),K,T,\beta,|I|\right)n^{\frac{-2\beta}{2\beta+1}}\]
with $C_1\left(q,\mathbb{E}\left(\|l_T\|_{I,\infty} | D\left(I,\nu\right)\right),K,T,\beta,|I|\right)$ and $C_2\left(q,\mathbb{E}\left(\|l_T\|_{I,\infty} | D\left(I,\nu\right)\right),K,T,\beta,|I|\right)$ constants depending on $q$, $\mathbb{E}\left(\|l_T\|_{I,\infty} | D\left(I,\nu\right)\right)$, $K$, $T$, $\beta$ and $|I|$.
Finally,
\[
\underset{n \to \infty}{\lim \sup}\; {\bf R}\left(\hat{q}_n, \Lambda_{\rho,\beta}, n^{\frac{-\beta}{2\beta+1}}\right) < \infty.
\]
In the case where $l > m$, we apply the Taylor expansion formula of the bias up to the order $m$ and we find that the bias is of oder $n^{2m}$ as $q^{(m)}$ is bounded. The convergence rate is then of order $n^{\frac{-m}{2m+1}}$.

\subsection{Proof of Proposition \ref{parametricconvergence}}

For simplicity, the proof is done in the case where $\Theta \subset \mathbb{R}$, that is when $d = 1$. The proof is similar for any $d \geq 1$. In the following, let us use the following notation:
\[\tilde{K}\left(x,h,z\right) = w\left(x,h,z\right)K\left(z\right), \text{ for } x \in I,\; h > 0,\; z \in \mathbb{R}.\]
As $X$ is an ancillary statistic, one can work as if $X$ was deterministic. We also work on the event $D\left(I,\nu\right)$.

\paragraph{{\bf Convergence of $\hat{\theta}_n$}} First, let us study the convergence of $M_n\left(\theta\right)$ for $\theta \in \Theta$. We can write $M_n\left(\theta\right)$ as the sum of 
\begin{equation} \label{martingale} \frac{2}{n^2 h_n^2} \int_0^T \int_0^{s^-} \int_I \tilde{K}\left(x,h_n,\frac{X_s-x}{h_n}\right)\tilde{K}\left(x,h_n,\frac{X_u-x}{h_n}\right){\bf 1}_{X_s \in I}{\bf 1}_{X_u \in I} dx dM_u dM_s,
\end{equation}
\begin{equation} \label{meanpart}
\frac{1}{h_n^2}\int_I \left( \int_0^T \tilde{K}\left(x,h_n,\frac{X_s-x}{h_n}\right){\bf 1}_{X_s \in I}\left(q\left(X_s\right)-g_{\theta}\left(X_s\right)\right)ds\right)^2dx   
\end{equation}
and 
\begin{equation} \label{crosstermcontrast}
\frac{2}{n h_n^2}\int_I \int_0^T \tilde{K}\left(x,h_n,\frac{X_s-x}{h_n}\right){\bf 1}_{X_s \in I} dM_s \int_0^T \tilde{K}\left(x,h_n,\frac{X_s-x}{h_n}\right) \left(q\left(X_s\right)-g_{\theta}\left(X_s\right)\right) {\bf 1}_{X_s \in I} ds dx.    
\end{equation}
The first term \eqref{martingale} has expectation $0$ and variance equal to 
\[\frac{4}{n^4 h_n^4} \int_0^T \int_0^{s^-} \left(\int_I \tilde{K}\left(x,h_n,\frac{X_s-x}{h_n}\right)\tilde{K}\left(x,h_n,\frac{X_u-x}{h_n}\right){\bf 1}_{X_s \in I}{\bf 1}_{X_u \in I} dx \right)^2 n^2 q\left(X_u\right)q\left(X_s\right)du ds\]
which is equal to 
\begin{equation} \label{variancedoubleintegral}
\frac{2}{n^2 h_n^4} \int_0^T \int_0^{T} \left(\int_I \tilde{K}\left(x,h_n,\frac{X_s-x}{h_n}\right)\tilde{K}\left(x,h_n,\frac{X_u-x}{h_n}\right){\bf 1}_{X_s \in I}{\bf 1}_{X_u \in I} dx\right)^2 q\left(X_u\right)q\left(X_s\right)du ds.
\end{equation}
Using the occupation time formula, \eqref{variancedoubleintegral} is equal to 
\[\frac{2}{n^2 h_n^4} \int_I \int_I  \left(\int_I \tilde{K}\left(x,h_n,\frac{r-x}{h_n}\right)\tilde{K}\left(x,h_n,\frac{y-x}{h_n}\right)dx\right)^2  q\left(r\right) l_T^r q\left(y\right)l_T^y dr dy.\]
By writing $I = \left[\underbar{I},\bar{I}\right]$, \eqref{variancedoubleintegral} is equal to 
\[\frac{2}{h_n n^2} \hspace{-0.2em}\int_I  \hspace{-0.2em}\int_{\frac{\underline{I}-y}{h_n}}^{\frac{\bar{I}-y}{h_n}}  \hspace{-0.5em} \left(\int_{\frac{y-\bar{I}}{h_n}}^{\frac{y-\underline{I}}{h_n}}\tilde{K}\left(y-uh_n,h_n,u\right)\tilde{K}\left(y-uh_n,h_n,u+p\right) du\right)^2 \hspace{-0.2em} q\left(y+h_n p\right) l_T^{y+h_n p} q\left(y\right)l_T^y dp dy.\]
Using the continuity properties of $x \mapsto l_T^x$ from Proposition \ref{propertiesX}, several times the dominated convergence theorem and the fact that $\underset{x \in I}{\inf} l_T^x \geq \nu > 0$ and $\|l_T\|_{I,\infty} < \infty$, we find that 
\[\tilde{K}\left(x,h_n,z\right) \underset{n \to \infty}{\longrightarrow} \left(l_T^x\right)^{-1}w\left(z\right) K\left(z\right), \text{ for }x\in I, \; z \in \mathbb{R}\]
with 
\[w\left(u\right) = U^T\left(0\right)\left(\int_{\mathbb{R}} U\left(z\right)U^T\left(z\right)K\left(z\right)dz\right)^{-1}U\left(u\right)\]
and that \eqref{variancedoubleintegral} is equivalent to 
\[\frac{2}{h_n n^2 } \int_{\mathbb{R}} \left(\int_{\mathbb{R}} w\left(u\right) w\left(u+p\right)K\left(u\right) K\left(u+p\right) du\right)^2 dp \int_I \frac{\left(q\left(y\right)\right)^2}{\left(l_T^y\right)^2} dy.\]
Thus, \eqref{martingale} is of order $O_p\left(\frac{1}{n \sqrt{h_n}}\right)$ and goes to 0 in probability. Concerning the second term \eqref{meanpart}, using the dominated convergence theorem again, we find that it converges to
\[\int_{\mathbb{R}} w\left(u\right)K\left(u\right)du \int_I \left(q\left(x\right)-g_\theta\left(x\right)\right)^2dx\]
and 
\[\int_{\mathbb{R}} w\left(u\right)K\left(u\right)du = 1\]
as it is the limit of $\int_{0}^T w\left(x,h_n,\frac{X_s-x}{h_n}\right)K_{h_n}\left(\frac{X_s-x}{h_n}\right)ds$ which is equal to 1 according to Proposition \ref{polynom}. The last term \eqref{crosstermcontrast} has mean 0 and variance equal to $O\left(\frac{1}{n}\right)$ and goes to 0 when $n \to \infty$.
Thus, 
\[M_n\left(\theta\right) \overset{p}{\rightarrow} M\left(\theta\right) = \| q - g_{\theta}\|_I^2.\]
Under Assumption \ref{assumpparametric}, we easily show that 
\[| M_n\left(\theta_1\right) - M\left(\theta_1\right) - \left(M_n\left(\theta_2\right) - M\left(\theta_2\right)\right) | \lesssim |\theta_1 - \theta_2|.\]
According to Kolmogorov continuity criterion, the convergence is then uniform in $\theta$. Furthermore, under $H_0$, the minimum of $M\left(\theta\right)$ is achieved for $\theta = \theta_0$ which is the unique minimum under Assumption \ref{assumpparametric}. We then have, under $H_0$, 
\[\hat{\theta}_n \overset{p}{\rightarrow} \theta_0.\]
Using Taylor's formula, on the event $\{M_n^{\left(2\right)}\left(\theta_0\right) + \frac{\left(\hat{\theta}_n-\theta_0\right)}{2}M^{\left(3\right)}\left(\tilde{\theta}_n\right) \neq 0\}$,
\[
\sqrt{n}\left(\hat{\theta}_n-\theta_0\right) = \frac{-\sqrt{n} M_n'\left(\theta_0\right)}{M_n^{\left(2\right)}\left(\theta_0\right) + \frac{\left(\hat{\theta}_n-\theta_0\right)}{2}M^{\left(3\right)}\left(\tilde{\theta}_n\right)}
\]
with $\tilde{\theta}_n$ between $\theta_0$ and $\hat{\theta}_n$. $M'_n\left(\theta\right)$ is equal to 
\[-\frac{2}{n h^2_n}  \int_0^T f\left(h_n,\theta,X_s\right) \left(dN_s-ng_{\theta}\left(X_s\right)ds\right) \]
with 
\[f\left(h_n,\theta,X_s\right) = \int_I \int_0^T \tilde{K}\left(x,h_n,\frac{X_u-x}{h_n}\right) {\bf 1}_{X_u \in I} \partial_{\theta} g_\theta\left(X_u\right) du \tilde{K}\left(x,h_n,\frac{X_s-x}{h_n}\right){\bf 1}_{X_s \in I} dx,\]
$M_n^{\left(2\right)}\left(\theta\right)$ to
\[ \begin{split}
& -\frac{2}{n h_n^2} \int_0^T\partial_{\theta} f\left(h_n,\theta,X_s\right)\left(dN_s-ng_{\theta}\left(X_s\right)ds\right) \\
&+ \frac{2}{h_n^2} \int_I\left(\int_0^T \tilde{K}\left(x,h_n,\frac{X_s-x}{h_n}\right) {\bf 1}_{X_s \in I} \partial_{\theta} g_\theta\left(X_s\right)  ds\right)^2dx
\end{split}
\]
and $M_n^{\left(3\right)}\left(\theta\right)$ to
\[ \begin{split}
& -\frac{2}{n h_n^2} \int_0^T \partial_{2,\theta} f\left(h_n,\theta,X_s\right)  \left(dN_s-ng_{\theta}\left(X_s\right)ds\right)\\
&+\frac{6}{h_n^2} \int_I \int_0^T \tilde{K}\left(x,h_n,\frac{X_s-x}{h_n}\right){\bf 1}_{X_s \in I} \partial_{\theta} g_\theta\left(X_s\right)  ds \int_0^T \tilde{K}\left(x,h_n,\frac{X_s-x}{h_n}\right) {\bf 1}_{X_s \in I} \partial_{2,\theta} g_\theta\left(X_s\right)  ds dx.
\end{split}
\]
Under $H_0$, $M_s = N_s - n\int_{0}^s g_{\theta_0}\left(X_s\right)ds$. Thus, $\sqrt{n} M'_n\left(\theta_0\right)$ has expectation $0$ and variance 
\[\frac{4}{h_n^4} \int_0^T \left(\int_I  \int_0^T \tilde{K}\left(x,h_n,\frac{X_s-x}{h_n}\right)\partial_\theta g_{\theta_0} \left(X_s\right) {\bf 1}_{X_s \in I}ds \tilde{K}\left(x,h_n,\frac{X_u-x}{h_n}\right) {\bf 1}_{X_u \in I} \right)^2 g_{\theta_0}\left(X_u\right) du\]
that converges to 
\[4\int_I \frac{\left(\partial_{\theta} g_{\theta_0}\left(u\right)\right)^2 g_{\theta_0}\left(u\right)}{l_T^u}du.\]
We can easily show that 
\[\frac{8 \left(\sqrt{n}\right)^3}{n^3 h_n^6} \int_0^T\left(f\left(h_n,\theta,X_s\right)\right)^3 n g_{\theta_0}\left(X_u\right) du = O\left(\frac{1}{\sqrt{n}}\right)\]
and converges to 0. Thus, according to \cite[Theorem 3]{peccati08}, we have the following central limit theorem
\begin{equation} \label{convergenceM1}
-\sqrt{n}M_n'\left(\theta_0\right) \overset{\mathcal{L}}{\rightarrow} \sqrt{4 \int_I \frac{\left(\partial_{\theta} g_{\theta_0}\left(u\right)\right)^2 g_{\theta_0}\left(u\right)}{l_T^u}du } \mathcal{N}\left(0,1\right).
\end{equation}
Under $H_0$, the first term of $M_n^{\left(2\right)}\left(\theta_0\right)$ 
\begin{equation} \label{M2term1}
 -\frac{2}{n h_n^2}  \int_0^T \partial_{\theta} f\left(h_n,\theta,X_s\right)  dM_s
 \end{equation}
has mean $0$ and variance 
\[ \frac{4}{n^2 h_n^4} \int_0^T \left(\partial_{\theta} f\left(h_n,\theta,X_s\right)\right)^2 n g_{\theta_0}\left(X_s\right)ds = O\left(\frac{1}{n}\right).\]
Then, \eqref{M2term1} converges to 0 in probability. The second term of $M^{\left(2\right)}\left(\theta_0\right)$ 
\[
 \frac{2}{h_n^2} \int_I \left(\int_0^T \tilde{K}\left(x,h_n,\frac{X_s-x}{h_n}\right){\bf 1}_{X_s \in I} \partial_{\theta} g_\theta\left(X_s\right)  ds\right)^2dx
\]
converges in probability to 
\[   
2 \int_I \left(\partial_{\theta} g_{\theta_0}\left(x\right)\right)^2 dx.
\]
Then, 
\begin{equation} \label{convergenceM2}
M_n^{\left(2\right)}\left(\theta_0\right) \overset{p}{\rightarrow}2\int_I \left(\partial_{\theta} g_{\theta_0}\left(u\right)\right)^2 du.
\end{equation}
The absolute value of $M^{\left(3\right)}\left(\tilde{\theta}_n\right)$ has mean equal to $O\left(1\right)$. As $\hat{\theta}_n - \theta \overset{p}{\rightarrow} 0$,
\begin{equation} \label{convergenceM3}
\left(\hat{\theta}_n-\theta\right)M_n^{\left(3\right)}\left(\tilde{\theta}_n\right) \overset{p}{\rightarrow} 0.
\end{equation}
Using \eqref{convergenceM1}, \eqref{convergenceM2} and \eqref{convergenceM3},
\[\sqrt{n}\left(\hat{\theta}_n-\theta_0\right)  \overset{\mathcal{L}}{\rightarrow} \frac{ \sqrt{ \int_I \frac{\left(\partial_{\theta} g_{\theta_0}\left(u\right)\right)^2 g_{\theta_0}\left(u\right)}{l_T^u}du }}{\int_I \left(\partial_{\theta} g_{\theta_0}\left(u\right)\right)^2 du}\mathcal{N}\left(0,1\right),\]
achieving the proof of (i).

\paragraph{{\bf Convergence of $M_n\left(\hat{\theta}_n\right)$ under $H_0$.}} $M_n\left(\hat{\theta}_n\right)$ is the sum of \eqref{martingale}, \eqref{meanpart} and \eqref{crosstermcontrast}, replacing $\theta$ by $\hat{\theta}_n$. Concerning \eqref{martingale}, we have 
\[\frac{16}{n^4 h_n^6} \int_0^T \int_0^{s^-} \left(\int_I \tilde{K}\left(x,h_n,\frac{X_s-x}{h_n}\right)\tilde{K}\left(x,h_n,\frac{X_u-x}{h_n}\right){\bf 1}_{X_s \in I}{\bf 1}_{X_u \in I} dx\right)^4 n^2  q\left(X_u\right)q\left(X_s\right)du ds
\]
which is equivalent to 
\[\frac{8}{h_n n^2} \int_{\mathbb{R}} \left(\int_{\mathbb{R}} w\left(u\right)w\left(u+p\right) K\left(u\right) K\left(u+p\right) du\right)^4 dp \int_I \frac{\left(q\left(y\right)\right)^2}{\left(l_T^y\right)^6} dy\]
and converges to 0. According to \cite[Theorem 3]{peccati08}, 
\[ \frac{2\sqrt{h_n} n}{n^2 h_n^2}  \int_0^T \int_0^{s^-} \int_I \tilde{K}\left(x,h_n,\frac{X_s-x}{h_n}\right)\tilde{K}\left(x,h_n,\frac{X_u-x}{h_n}\right) dx {\bf 1}_{X_s \in I}{\bf 1}_{X_u \in I} dM_u dM_s\]
converges in law to
\[ \sqrt{2 \int_{\mathbb{R}} \left(\int_{\mathbb{R}} w\left(u\right)w\left(u+p\right)K\left(u\right) K\left(u+p\right) du\right)^2 dp \int_I \frac{\left(q\left(y\right)\right)^2}{\left(l_T^y\right)^2} dy} \mathcal{N}\left(0,1\right).\]
The term \eqref{meanpart} is of the same order than $\left(\hat{\theta}_n-\theta_0\right)^2$ which is $O_p\left(\frac{1}{n}\right)$ and 
\[ \frac{\sqrt{h_n} n}{h_n^2}  \int_I \left( \int_0^T \tilde{K}\left(x,h_n,\frac{X_s-x}{h_n}\right){\bf 1}_{X_s \in I}\left(q\left(X_s\right)-g_{\hat{\theta}_n}\left(X_s\right)\right)ds\right)^2dx  = O_p\left( \sqrt{h_n}\right) \overset{p}{\rightarrow} 0.\]
The last term corresponding to \eqref{crosstermcontrast} is bounded by 
\[\frac{2}{nh_n^2}M_I |\theta_0-\hat{\theta}_n| \int_I \int_{0}^T \tilde{K}\left(x,h_n,\frac{X_s-x}{h_n}\right) {\bf 1}_{X_s \in I} dM_s \int_0^T  \tilde{K}\left(x,h_n,\frac{X_s-x}{h_n}\right)  {\bf 1}_{X_s \in I} ds\]
equal to $O_p\left(\frac{1}{n}\right)$.
Thus, 
\[ \frac{2\sqrt{h_n}n}{n h_n^2}\int_I \int_0^T \tilde{K}\left(x,h_n,\frac{X_s-x}{h_n}\right) {\bf 1}_{X_s \in I} dM_s \int_0^T  \tilde{K}\left(x,h_n,\frac{X_s-x}{h_n}\right) \left(q\left(X_s\right)-g_{\hat{\theta}_n}\left(X_s\right)\right) {\bf 1}_{X_s \in I} ds\]
is equal to $O_p\left(\sqrt{h_n}\right)$ and converges to 0 in probability. Finally, 
\[n\sqrt{h_n}M_n\left(\hat{\theta}_n\right)  \overset{\mathcal{L}}{\rightarrow}   \sqrt{2 \int_{\mathbb{R}} \left(\int_{\mathbb{R}}w\left(u\right)w\left(u+p\right) K\left(u\right) K\left(u+p\right) du\right)^2 dp \int_I \frac{\left(q\left(y\right)\right)^2}{\left(l_T^y\right)^2} dy} \mathcal{N}\left(0,1\right).\]
 	
\paragraph{{\bf Convergence of $M_n\left(\hat{\theta}_n\right)$ under $H_1$.}} As the convergence of $M_n\left(\theta\right)$ is uniform in $\theta$, $M_n\left(\hat{\theta}_n\right)$ converges to $\underset{\theta \in \Theta}{\inf} \|q-g_{\theta}\|_I^2$ which is strictly positive under $H_1$. Thus, under $H_1$,
\[|n \sqrt{h_n} M_n\left(\hat{\theta}_n\right)| \rightarrow \infty.\]

\vip

\noindent {\bf Acknowledgements.} I am grateful to Olivier F\'eron and Marc Hoffmann for helpful discussion and comments. This research is supported by the department OSIRIS (Optimization, SImulation, RIsk and Statistics for Energy Markets) of EDF in the context of a CIFRE contract and by FiME (Finance for Energy Markets) Research Initiative.

\bibliographystyle{plain}
\bibliography{Biblio}

\begin{thebibliography}{10}

\bibitem{aalen78}
Odd Aalen.
\newblock Nonparametric inference for a family of counting processes.
\newblock {\em The Annals of Statistics}, pages 701--726, 1978.

\bibitem{ait96}
Yacine Ait-Sahalia.
\newblock Testing continuous-time models of the spot interest rate.
\newblock {\em Review of Financial studies}, 9(2):385--426, 1996.

\bibitem{ait15}
Yacine A{\"\i}t-Sahalia, Julio Cacho-Diaz, and Roger~JA Laeven.
\newblock Modeling financial contagion using mutually exciting jump processes.
\newblock {\em Journal of Financial Economics}, 117(3):585--606, 2015.

\bibitem{barlow82}
Martin~T Barlow and Marc Yor.
\newblock Semi-martingale inequalities via the garsia-rodemich-rumsey lemma,
  and applications to local times.
\newblock {\em Journal of functional Analysis}, 49(2):198--229, 1982.

\bibitem{benth12}
Fred~Espen Benth, R{\"u}diger Kiesel, and Anna Nazarova.
\newblock A critical empirical study of three electricity spot price models.
\newblock {\em Energy Economics}, 34(5):1589--1616, 2012.

\bibitem{benth15}
Fred~Espen Benth, Nina Lange, and Tor~Age Myklebust.
\newblock Pricing and hedging quanto options in energy markets.
\newblock {\em Journal of Energy Markets}, 2015.

\bibitem{benth11}
Fred~Espen Benth and J{\=u}rat{\.e} {\v{S}}altyt{\.e}~Benth.
\newblock Weather derivatives and stochastic modelling of temperature.
\newblock {\em International Journal of Stochastic Analysis}, 2011, 2011.

\bibitem{brooks91}
Maria~Mori Brooks and J~Stephen Marron.
\newblock Asymptotic optimality of the least-squares cross-validation bandwidth
  for kernel estimates of intensity functions.
\newblock {\em Stochastic Processes and their Applications}, 38(1):157--165,
  1991.

\bibitem{cartea05}
Alvaro Cartea and Marcelo~G Figueroa.
\newblock Pricing in electricity markets: a mean reverting jump diffusion model
  with seasonality.
\newblock {\em Applied Mathematical Finance}, 12(4):313--335, 2005.

\bibitem{castellan00}
G~Castellan and F~Letu{\'e}.
\newblock Estimation of the cox regression function via model selection.
\newblock {\em F. Letu{\'e}'s PhD Thesis}, 2000.

\bibitem{chen11}
Feng Chen, Paul~SF Yip, and KF~Lam.
\newblock On the local polynomial estimators of the counting process intensity
  function and its derivatives.
\newblock {\em Scandinavian Journal of Statistics}, 38(4):631--649, 2011.

\bibitem{comte11}
Fabienne Comte, St{\'e}phane Ga{\"\i}ffas, and Agathe Guilloux.
\newblock Adaptive estimation of the conditional intensity of marker-dependent
  counting processes.
\newblock In {\em Annales de l'institut Henri Poincar{\'e} (B)}, volume~47,
  pages 1171--1196, 2011.

\bibitem{cox92}
David~R Cox.
\newblock Regression models and life-tables.
\newblock In {\em Breakthroughs in statistics}, pages 527--541. Springer, 1992.

\bibitem{delattre13}
Sylvain Delattre, Christian~Y Robert, and Mathieu Rosenbaum.
\newblock Estimating the efficient price from the order flow: a brownian cox
  process approach.
\newblock {\em Stochastic Processes and their Applications}, 123(7):2603--2619,
  2013.

\bibitem{deschatre18}
Thomas Deschatre, Olivier F\'eron, and Marc Hoffmann.
\newblock Estimating fast mean-reverting jumps in electricity market models.
\newblock {\em arXiv preprint arXiv:1803.03803}, 2018.

\bibitem{diggle85}
Peter Diggle.
\newblock A kernel method for smoothing point process data.
\newblock {\em Applied statistics}, pages 138--147, 1985.

\bibitem{goldenshluger11}
Alexander Goldenshluger and Oleg Lepski.
\newblock Bandwidth selection in kernel density estimation: oracle inequalities
  and adaptive minimax optimality.
\newblock {\em The Annals of Statistics}, pages 1608--1632, 2011.

\bibitem{hoffmann01}
Marc Hoffmann.
\newblock On estimating the diffusion coefficient: parametric versus
  nonparametric.
\newblock In {\em Annales de l'IHP Probabilit{\'e}s et statistiques},
  volume~37, pages 339--372, 2001.

\bibitem{houdre03}
Christian Houdr{\'e} and Patricia Reynaud-Bouret.
\newblock Exponential inequalities, with constants, for u-statistics of order
  two.
\newblock In {\em Stochastic inequalities and applications}, pages 55--69.
  Springer, 2003.

\bibitem{karr91}
Alan Karr.
\newblock {\em Point processes and their statistical inference}, volume~7.
\newblock CRC press, 1991.

\bibitem{lacour16}
Claire Lacour, Pascal Massart, and Vincent Rivoirard.
\newblock Estimator selection: a new method with applications to kernel density
  estimation.
\newblock {\em Sankhya A}, 79(2):298--335, 2017.

\bibitem{lemler16}
Sarah Lemler et~al.
\newblock Oracle inequalities for the lasso in the high-dimensional aalen
  multiplicative intensity model.
\newblock In {\em Annales de l'Institut Henri Poincar{\'e}, Probabilit{\'e}s et
  Statistiques}, volume~52, pages 981--1008. Institut Henri Poincar{\'e}, 2016.

\bibitem{lerasle16}
Matthieu Lerasle, Nelo~Molter Magalh{\~a}es, and Patricia Reynaud-Bouret.
\newblock Optimal kernel selection for density estimation.
\newblock In {\em High Dimensional Probability VII}, pages 425--460. Springer,
  2016.

\bibitem{meyer08}
Thilo Meyer-Brandis and Peter Tankov.
\newblock Multi-factor jump-diffusion models of electricity prices.
\newblock {\em International Journal of Theoretical and Applied Finance},
  11(05):503--528, 2008.

\bibitem{murphy91}
Susan~Allbritton Murphy and Pranab~Kumar Sen.
\newblock Time-dependent coefficients in a cox-type regression model.
\newblock {\em Stochastic Processes and their Applications}, 39(1):153--180,
  1991.

\bibitem{nielsen95}
Jens~P Nielsen and Oliver~B Linton.
\newblock Kernel estimation in a nonparametric marker dependent hazard model.
\newblock {\em The Annals of Statistics}, pages 1735--1748, 1995.

\bibitem{ogata86}
Yosihiko Ogata and Koichi Katsura.
\newblock Point-process models with linearly parametrized intensity for
  application to earthquake data.
\newblock {\em Journal of applied probability}, 23(A):291--310, 1986.

\bibitem{osullivan93}
Finbarr O'Sullivan.
\newblock Nonparametric estimation in the cox model.
\newblock {\em The Annals of Statistics}, pages 124--145, 1993.

\bibitem{peccati08}
Giovanni Peccati, Murad~S Taqqu, et~al.
\newblock Central limit theorems for double poisson integrals.
\newblock {\em Bernoulli}, 14(3):791--821, 2008.

\bibitem{revuz13}
Daniel Revuz and Marc Yor.
\newblock {\em Continuous martingales and Brownian motion}, volume 293.
\newblock Springer Science \& Business Media, 2013.

\bibitem{reynaud03}
Patricia Reynaud-Bouret.
\newblock Adaptive estimation of the intensity of inhomogeneous poisson
  processes via concentration inequalities.
\newblock {\em Probability Theory and Related Fields}, 126(1):103--153, 2003.

\bibitem{reynaud14}
Patricia Reynaud-Bouret.
\newblock Concentration inequalities, counting processes and adaptive
  statistics.
\newblock In {\em ESAIM: Proceedings}, volume~44, pages 79--98. EDP Sciences,
  2014.

\bibitem{tankov03}
Peter Tankov.
\newblock {\em Financial modelling with jump processes}, volume~2.
\newblock CRC press, 2003.

\bibitem{truccolo05}
Wilson Truccolo, Uri~T Eden, Matthew~R Fellows, John~P Donoghue, and Emery~N
  Brown.
\newblock A point process framework for relating neural spiking activity to
  spiking history, neural ensemble, and extrinsic covariate effects.
\newblock {\em Journal of neurophysiology}, 93(2):1074--1089, 2005.

\bibitem{tsybakov09}
Alexandre~B Tsybakov.
\newblock Introduction to nonparametric estimation, 2009.

\bibitem{utikal93}
Klaus~J Utikal.
\newblock Nonparametric inference for a doubly stochastic poisson process.
\newblock {\em Stochastic processes and their applications}, 45(2):331--349,
  1993.

\bibitem{varet17}
Suzanne Varet, Claire Lacour, Pascal Massart, and Vincent Rivoirard.
\newblock Performances num{\'e}riques du choix de fen{\^e}tre par pco lors de
  l'estimation de densit{\'e}s multivari{\'e}es par m{\'e}thode {\`a} noyau.
\newblock 2017.

\bibitem{yen13}
Ju-Yi Yen and Marc Yor.
\newblock Local times and excursion theory for brownian motion a tale of wiener
  and ito measures preface, 2013.

\bibitem{zhang10}
Tingting Zhang and SC~Kou.
\newblock Nonparametric inference of doubly stochastic poisson process data via
  the kernel method.
\newblock {\em The annals of applied statistics}, 4(4):1913, 2010.

\end{thebibliography}
\end{document}